\documentclass{article}
\usepackage{epsfig}
\usepackage{color}
\usepackage{amsmath,amssymb,amscd,amsthm,bm} 
\usepackage{mathdots} 
 
\title{Continued fractions and Hankel determinants from hyperelliptic curves}

\author{Andrew N.W. Hone\thanks{On leave at 
School of Mathematics \&  
Statistics, University of New South Wales, 
NSW 2052, Australia.}~\\
School of Mathematics, 
Statistics \& Actuarial Science~\\ 
University of
Kent~\\ Canterbury CT2 7NF, UK.
}

\newcommand{\beq}{\begin{equation}}  
\newcommand{\eeq}{\end{equation}}  
\newcommand{\bea}{\begin{eqnarray}} 
\newcommand{\eea}{\end{eqnarray}}   
\newcommand{\bear}{\begin{array}}  
\newcommand{\eear}{\end{array}}

\newtheorem{thm}{Theorem}[section] 

\newtheorem{propn}[thm]{Proposition} 

\newtheorem{conje}[thm]{Conjecture}

\newenvironment{prf}{\trivlist \item [\hskip 
\labelsep {\bf Proof:}]\ignorespaces}{\qed \endtrivlist}

\theoremstyle{definition}

\newtheorem{exa}{Example}[section] %
\newtheorem{remark}{Remark}[section]

\newcommand{\Q}{{\mathbb Q}}
\newcommand{\Z}{{\mathbb Z}}
\newcommand{\C}{{\mathbb C}}
\newcommand{\R}{{\mathbb R}}

\newcommand{\rR}{\mathrm{R}}
\newcommand{\rQ}{\mathrm{Q}}
\newcommand{\rP}{\mathrm{P}}

\newcommand{\rd}{\mathrm{d}}
\newcommand{\ri}{\mathrm{i}}
\newcommand{\rg}{\mathrm{g}}
\newcommand\la{{\lambda}}

\newcommand\al{{\alpha}}
\newcommand\be{{\beta}}
\newcommand\ze{{\zeta}}
\newcommand\gam{{\gamma}}

\newcommand\om{{\omega}}

\newcommand\si{{\sigma}}

\def\lc{\left\lfloor}   
\def\rc{\right\rfloor}

\begin{document} 

\maketitle

\normalsize
 \begin{abstract} 
\noindent 
Following van der Poorten, 
we consider a family of nonlinear maps which are generated from the continued fraction 
expansion of a function on a hyperelliptic curve of genus $\mathrm{g}$. Using the connection 
with the classical theory of J-fractions and orthogonal polynomials, we 
show that in the simplest case $\mathrm{g}=1$ this provides a straightforward 
derivation of Hankel determinant formulae for the terms of a general Somos-4 sequence, 
which were found in a particular form by Chang, Hu and Xin, 
We extend these formulae to the higher genus case, and prove that generic Hankel 
determinants  in genus two satisfy a Somos-8 relation. Moreover, for all $\mathrm{g}$ we 
show that the iteration for the continued fraction expansion is equivalent to a discrete Lax pair 
with a natural Poisson structure, and the associated nonlinear map is a 
discrete integrable system.  \\ \\
\textit{Dedicated to the memory of Jon Nimmo.}\\ \\
\noindent {\bf Keywords: J-fraction, hyperelliptic curve, Hankel determininant, Somos sequence.} 
\end{abstract} 
\section{Introduction} 

\setcounter{equation}{0}

The Somos-4 recurrence is given by 
\beq \label{s4recu} 
\tau_{n+4}\tau_n =\alpha \, \tau_{n+3}\tau_{n+1} + \beta\, \tau_{n+2}^2 . 
\eeq 
The surprising observation of Somos was that when $\al=\be=1$ and the 
four initial values $\tau_0,\tau_1,\tau_2,\tau_3$ are all 1, the recurrence (\ref{s4recu}) generates a sequence of 
integers \cite{somos}, beginning with 
\beq\label{s4seq} 
1,1,1,1,2,3,7,23,59,314,1529,8209,83313, 620297, \ldots . 
\eeq 
A proof of this fact was eventually published \cite{malouf}, 
but a better understanding of the mechanism by which such 
rational  recurrences can yield integer sequences came 
from the observation that (\ref{s4recu}) exhibits the Laurent property \cite{gale}: 
the iterates are Laurent polynomials in the initial values with integer coefficients, 
that is to say 
$$ 
\tau_n\in \Z [\tau_0^{\pm 1},  \tau_1^{\pm 1},\tau_2^{\pm 1},\tau_3^{\pm 1},\al,\be] 
$$ 
for all $n$, which makes it obvious why (\ref{s4seq}) consists 
entirely of integers.  The Laurent phenomenon \cite{fz} eventually appeared 
as a key property of the distinguished generators (cluster variables) in Fomin and Zelevinsky's 
cluster algebras \cite{fz1}, which are constructed by a recursive process called mutation, 
and Fordy and Marsh showed how the Somos-4 recurrence and 
various higher order analogues arise from cluster mutations 
starting from quivers of a particular type \cite{fordy_marsh}.  

Cluster algebras fit within a broader setting 
of Laurent phenomenon algebras \cite{lp}, leading to a wide variety of nonlinear recurrences that 
exhibit the Laurent property \cite{alman}; 
and it is possible to reverse engineer a rational recurrence to generate 
an integer sequence \cite{zeil}. However, Somos-4 sequences have some 
very special features which are a consequence of the fact that each such 
sequence is associated with a sequence of points $P_0+nP$ on an elliptic curve $E$, 
and this leads to 
an analytic formula for the terms of the sequence. The 
following result was proved in \cite{honeblms}. 

\begin{thm}\label{s4sigma} 
The terms of a Somos-4 sequence, generated by (\ref{s4recu}) from 
four non-zero initial values $\tau_0,\tau_1,\tau_2,\tau_3$, non-zero 
$\al$ 
and arbitrary $\be\in\C$, are given by 
\beq\label{s4sig} 
\tau_n=\hat{a}\hat{b}^n\, \frac{\si (z_0+nz)}{\si (z)^{n^2}}, 
\eeq  
where $\si (z)=\si (z;g_2,g_3)$ is the 
Weierstrass sigma function associated with the elliptic curve 
$E: \, y^2=4x^3-g_2x-g_3$ over $\C$ with period lattice $\Lambda$, with  
$z_0=\int_\infty^{P_0}\frac{\rd x}{y}$,  $z=\int_\infty^{P}\frac{\rd x}{y}\in\C\bmod\Lambda$ 
corresponding to points $P_0,P\in E$, and $\hat{a},\hat{b}$ are certain non-zero constants.  
\end{thm} 

The formula (\ref{s4sig}) also makes 
sense in the degenerate case when the discriminant 
$g_2^3-27g_3^2=0$. 
Although the above result is formulated over the complex numbers, its algebraic 
content - associating a solution of (\ref{s4recu}) with a sequence of points 
on an elliptic curve - is valid in any field over which the initial values and 
coefficients $\al,\be$ are defined (up to appropriate adjustments in characteristic 2 or 3); this was described independently 
by Swart \cite{swart}, and also, in terms of a quartic model for $E$, by 
van der Poorten \cite{vdp}. 
This underlying algebraic structure has many consequences, 
including the existence of higher order relations between the terms 
\cite{hones5, magic, swartvdp}, and more refined versions of the Laurent 
property which produce large families of integer sequences \cite{swahon}. 
From this point of view, Somos-4 sequences are natural extensions 
of Ward's elliptic divisibility sequences \cite{ward}, which correspond to the 
special case $P_0=\infty$ (the identity element in the group law of $E$), 
and generalize the arithmetical properties of Fibonacci or Lucas sequences 
to a nonlinear setting \cite{recs}. Aside from  their intrinsic interest 
for certain problems of an arithmetical \cite{rob} or Diophantine nature \cite{heron}, Somos 
sequences and their higher order analogues appear in discrete integrable systems,  
underlying many integrable maps \cite{hkq}, especially via reductions 
of the discrete Hirota equation (bilinear discrete KP, also known as the octahedron recurrence) \cite{hkw} 
or Miwa's equation  (bilinear discrete BKP, or the cube recurrence) \cite{s6prym}. 
They also arise in solvable models of statistical mechanics and quantum field theory, 
such as the hard hexagon model, as mentioned in \cite{qrt}, or dimer models 
and quiver gauge theory \cite{eager}.

It was conjectured by Somos, and later proved by Xin \cite{xin}, that the terms of the 
sequence (\ref{s4seq}) have another explicit expression that is rather different from 
(\ref{s4sig}), being given by the Hankel determinants 
\beq\label{origh} 
D_n =\det (\tilde{s}_{i+j-2})_{i,j=1,\ldots, n}, 
\eeq 
where the entries $\tilde{s}_j$ are obtained from the function 
$\eta = \eta (x)$ satisfying the algebraic equation 
\beq\label{37a}
\eta-\eta^2 = x-x^3.  
\eeq  
To be precise, solving (\ref{37a}) for $\eta$ with a fixed choice of square root, one should take the function 
$\tilde{G} =\eta/x -1$, 
and expand it as 
\beq \label{gtilde}
\frac{\eta}{x}-1 = \sum_{j\geq 1} \tilde{s}_{j-1}x^j  = x+x^2+3x^3 +8x^4+23x^5+\cdots,  
\eeq 
which gives 
\beq \label{dhank}
D_0=D_1=1, \qquad D_2 =  \left|\begin{array}{cc} 1 & 1\\ 1 & 3 \end{array}\right|=2, \quad 
D_3 = \left|\begin{array}{ccc} 1 & 1 & 3 \\ 1 & 3 & 8 \\ 3 & 8 & 23 \end{array}\right|=3, 
\eeq 
and so on, where the matrix entries are generated by the recursion 
\beq\label{g1altrec} 
\tilde{s}_j=\tilde{\al}\tilde{s}_{j-1}+ \tilde{\be}\tilde{s}_{j-2}
+\tilde{\gam}\sum_{i=0}^{j-2}\tilde{s}_{i}\tilde{s}_{j-2-i}, \qquad j\geq 2, 
\eeq 
with $\tilde{\alpha}=2$, $\tilde{\be}=0$, $\tilde{\gam}=1$ and $\tilde{s}_0=\tilde{s}_1=1$. 
It was further conjectured by Barry \cite{barry} that the Hankel determinants $D_n$ formed from a particular 
family of sequences $(\tilde{s}_j)$, defined by the recursion (\ref{g1altrec}) with 
$\tilde{s}_0=1$, $\tilde{s}_1=\tilde{\al}$, satisfy the Somos-4 recurrence 
(\ref{s4recu}) with coefficients $\al = \tilde{\al}^2\tilde{\gam}^2$, 
$\be = \tilde{\gam}^2(\tilde{\be}+\tilde{\gam})^2- \tilde{\al}^2\tilde{\gam}^3$, 
and this was proved by Chang and Hu using identities for block Hankel determinants \cite{ch}. 
The latter  result does not overlap with that of Xin, since the conditions on 
the coefficients and initial conditions do not include 
the  original sequence (\ref{s4seq}). However, it was subsequently shown 
by Chang, Hu and Xin that, for any Somos-4 sequence with two adjacent initial values 
equal to 1, the terms with positive index $n$ are given by a Hankel determinant of 
the form (\ref{origh}), where the entries $\tilde{s}_j$ satisfy a recursion of the form 
(\ref{g1altrec}), for a suitable   choice of 
$\tilde{\alpha}, \tilde{\be},\tilde{\gam},\tilde{s}_0,\tilde{s}_1$ \cite{chx}. 

In this paper we start from van der Poorten's construction in \cite{vdp} for Somos-4, 
based on the continued fraction expansion of a function on a quartic curve of genus one, and the 
results of \cite{vdp2,vdp3}. 
In  the latter work, the continued fraction approach was extended to hyperelliptic curves of higher genus $\rg$, 
defined by a polynomial of  even degree $2\rg+2$,   but 
with partial success: a Somos-6 relation was obtained in genus two, but only in a special case. 
The continued fraction expansion of a hyperelliptic function of a certain type is described in 
section \ref{cfs}  (for other related results, and the connection 
with the geometry of the Jacobian of the curve, see \cite{ar,berry, bombcoh,grossetv,vdphyp}). 
Next,  in section \ref{laxnonl},  the recursion for the  continued fraction is reformulated as a discrete 
dynamical system defined by a matrix linear problem (a discrete Lax pair), and we state the 
first main result, Theorem \ref{int}, which says that this nonlinear dynamical system is integrable in the sense 
that it satisfies a discrete analogue of the Liouville-Arnold theorem from classical mechanics \cite{arnold}, 
having an invariant symplectic structure and a sufficient number of first integrals in involution  
with respect to the corresponding Poisson bracket 
(see \cite{bruschi, maeda, veselov} for the precise definition of a discrete integrable system). 
We present explicit details of the 
symplectic map and first integrals for the cases $\rg=1$ and $\rg=2$, 
while the complete proof for any $\rg$ is deferred until section \ref{poisson}. 

Section \ref{hankels} is concerned with the derivation of Hankel determinant formulae 
for the solutions of the nonlinear system, based on the classical theory of J-fractions and 
orthogonal polynomials, which are presented in a uniform fashion for any genus $\rg$. 
(Note that Hankel determinants and continued fractions have appeared in
the solution of many other integrable 
systems, particularly those of Toda type, and  Painlev\'e equations \cite{commonh, jkm, suris}, 
and there are  more recent results in the broader context of Pad\'{e} approximants  associated with isomonodromic 
deformations  \cite{manotsuda}.) 
Subsequently, in section \ref{somoscon}, we show how the Hankel 
formulae obtained generalize the results of 
Chang, Hu and Xin on Somos-4 sequences: even in the elliptic case $\rg=1$, the results 
are more general than \cite{chx}, since the Hankel determinants  depend on an 
additional free parameter, and the formulae from the previous section extend to negative 
indices $n$. In genus two we prove that, for  generic parameter values,  the corresponding Hankel 
determinants satisfy a Somos-8 relation  (Theorem \ref{g2h}), 
and indicate how van der Poorten's  Somos-6 recurrence   arises as a special case. 
We also present a precise conjecture that provides an analytic formula analogous to 
(\ref{s4sig})  for the Hankel determinants in genus two, 
and briefly explain how an appropriate higher genus analogue of this conjecture implies 
the existence of Somos recurrences for all values of $\rg$.     
 
In section \ref{poisson} we employ  a space of $2\times 2$ Lax matrices, related to those in section \ref{laxnonl} 
by a gauge transformation, which admits a natural Poisson structure, and construct a completely 
integrable system on this phase space, given by a set of commuting flows defined by
suitable Hamiltonian functions. 
We then show how the nonlinear map coming from 
the continued fraction arises from a Poisson map on this phase space, which preserves the   
same Hamiltonians  and Casimirs as the continuous system. The map we obtain is somewhat 
reminiscent of  the B\"{a}cklund transformation (BT) for the  even Mumford systems, 
introduced in \cite{kuzvan}, except that the entries of the Lax pair have a different degree structure, 
and (in contast with the BT, which is a multivalued correspondence) it is an 
explicit birational map. Some determinantal identities that directly yield the formulae 
for the coefficients of J-fractions are also presented in an appendix.

\section{Continued fractions for hyperelliptic functions} \label{cfs}

\setcounter{equation}{0}

Following van der Poorten \cite{vdp,vdp2,vdp3}, we consider 
a hyperelliptic curve defined by 
\beq\label{curve} 
{\cal C}: \qquad Y^2=F(X), 
\eeq 
where 
\beq\label{fpoly} 
F(X)= A(X)^2+4R(X)
\eeq 
for a pair of polynomials 
$$A(X)=X^{\rg+1}+\cdots, \qquad R(X)=u\, X^\rg+\cdots. 
$$ 
In addition to the affine points $(X,Y)\in\C^2$ satisfying (\ref{curve}), 
one can adjoin two points at infinity, $\infty_1$ and $\infty_2$, such that,  
in terms of a local parameter $t$ at each of these points,  
$X=1/t$ and $Y\sim \pm 1/t^{\rg+1}$ respectively. Thus, in the generic case 
that all the roots of $F(X)$ are distinct, one obtains 
a compact Riemann surface of genus $\rg$, also denoted ${\cal C}$.  
In the associated 
function field 
${\cal F}=\C (X,Y) / (\,Y^2=A(X)^2+4R(X)\,)$, 
we pick 
\beq\label{pick} Y_0=\frac{Y+P_0}{Q_0}\in{\cal F},
\eeq  
for certain polynomials $P_0$, $Q_0$ in $X$ of degrees 
$\rg+1$, $\rg$  respectively, taking the form  
\beq\label{pqform} 
P_0(X)=A(X)+2d_0\,  X^{\rg-1}+\cdots,  
\qquad
Q_0 (X)= u_0\,  X^\rg+\cdots, \quad u_0\neq 0,  
\eeq 
and we impose the additional requirement that 
\beq\label{reduced}
Q_0(X)|Y^2-P_0(X)^2=F(X)-P_0(X)^2. 
\eeq 

To compute the continued fraction of any  element  $Y_0\in{\cal F}$, we 
take its  expansion in the neighbourhood of the point $\infty_1\in{\cal C}$, given by a power series in  
 $X^{-1}$;   
this can be viewed as an element of $\C((X^{-1}))$. Then the continued 
fraction is 
\beq\label{cfrac} 
Y_0=
a_0+\cfrac{ 1}{ a_1 +\cfrac {1}{ a_2 +\cdots} } 
=\left \lfloor{Y_0}\right \rfloor + \mathrm{remainder}, 
\eeq
where, for any element of $\C((X^{-1}))$, the floor symbol denotes the polynomial part, 
and the remainder is a series in positive powers of $X^{-1}$. 
Thus, by iterating   the standard recursion 
\beq\label{rec}
Y_n =a_n +\frac{1}{Y_{n+1}}, \qquad a_n =\lc Y_n \rc 
\eeq 
for $n=0,1,2,\ldots$, one obtains the successive partial quotients $a_n(X)$ in the continued fraction 
(\ref{cfrac})  above, which are 
polynomials in $X$.

Now let us describe in detail the form of the continued fraction expansion for the 
particular type of function $Y_0$ specified by (\ref{pick}). Because the neighbourhood 
of $\infty_1$ is being considered, with $X\to\infty$, it follows that 
$Y\sim A(X)\sim X^{\rg+1}$, so $Y+P_0\sim 2A(X)$ and hence $Y_0\sim 2X/u_0$. Thus $a_0$ 
is linear in $X$. Moreover, by (\ref{reduced}), there is 
some   polynomial $Q_{-1}(X)$ of degree $\rg$  such that 
$Y^2-P_0^2 =Q_0Q_{-1}$, so 
$$ 
Y_0=a_0+\frac{1}{Y_1} = \frac{Y+P_0}{Q_0} = \frac{Q_{-1}}{Y-P_0}. 
$$  
If a polynomial $P_1(X)$ of degree $\rg+1$ is defined by 
$$ 
P_1=-P_0 +a_0Q_0, 
$$ 
then it follows from (\ref{reduced}) that $Q_0|Y^2-P_1^2$, so there 
is some polynomial $Q_1(X)$ of degree $\rg$ such that 
\beq\label{y1} 
Y_1= \frac{Y+P_1}{Q_1} = \frac{Q_{0}}{Y-P_1}, 
\eeq 
and 
also $Q_1|Y^2-P_1^2$. 
Then by induction, at each stage of the recursion 
(\ref{rec}) we find 
linear 
partial quotients 
\beq\label{parq} 
a_n(X)=
{ 2(X+v_n)}  /{u_n}, 
\eeq  
and  
\beq\label{eig}
Y_n=\frac{Y+P_n}{Q_n} = \frac{Q_{n-1}}{Y-P_n},
\eeq 
for a sequence of polynomials $P_n$, $Q_n$ of degrees $\rg+1$ and $\rg$ respectively, 
with 
\beq\label{PQform} 
P_n(X)=A(X)+2d_n\,  X^{\rg-1}+\cdots,  
\qquad
Q_n(X)= u_n\,  X^\rg+\cdots . 
\eeq  
Note from (\ref{parq}) that $u_n\neq 0$ is required for the recursion to make 
sense at each stage.   Moreover, 
at each stage, $Y_n$ has positive degree in $X$, and $\bar{Y}_n =(-Y+P_n)/Q_n$ (its  
image under the hyperelliptic involution) has negative degree; in the terminology 
of van der Poorten, $Y_n$ is \textit{reduced} \cite{vdp3}. 

Observe that in the equation (\ref{curve}) for ${\cal C}$ there is always the freedom to shift 
$X\to X+\,$const, which replaces $F(X)$ by another monic polynomial of the same degree. 
Henceforth we will exploit this freedom  in order to remove 
the coefficient at order $X^\rg$ in $F$, which means that 
\beq\label{aform} 
A(X) = X^{\rg+1}  + \sum_{j=0}^{\rg-1}k^{(j)} X^j, 
\eeq 
for some constants $k^{(j)}$. 
This choice is convenient because in the continued fraction expansion it means that 
\beq\label{qv} 
Q_n(X)=u_n\Big(X^\rg - v_n X^{\rg -1}+ O(X^{\rg-2})\Big).
\eeq 
In other words, modulo  factors of 2 and  $u_n$, the next to leading order term in $Q_n$ 
completely fixes the constant term in the partial quotient (\ref{parq}), and we will always assume that 
$A$ has the form (\ref{aform}) so that this is the case. (We will make another comment 
about this later, when we discuss Poisson brackets.) 
 
As we shall see in the next section, the recursion  (\ref{rec}) and the relations (\ref{eig}) 
together yield a set of coupled nonlinear recurrences for the coefficients appearing 
in $P_n$ and $Q_n$. For the time being we derive just one such relation, by considering 
the second equality in (\ref{eig}), which is equivalent to the identity 
\beq\label{dety} 
Y^2=P_n^2 +Q_nQ_{n-1}.
\eeq 
If we  make use of (\ref{curve}) together with 
(\ref{PQform}), and cancel $A^2$ from both sides, then we have 
\beq\label{specrel} 
\begin{array}{rcl} 
4R(X) & = & 2A(X)\, \Big(2d_nX^{\rg-1} +O(X^{\rg-2})\Big) +\Big(4d_n^2X^{2\rg-2}+ 
O(X^{2\rg-3})\Big) \\ 
&& +u_nu_{n-1} 
\Big(X^\rg + O(X^{\rg-1})\Big)\,  \Big(X^\rg + O(X^{\rg-1})\Big), 
\end{array} \eeq 
so the leading order term, at order $X^{2\rg}$, gives the formula 
\beq \label{ud} 
4d_n +u_nu_{n-1}=0.
\eeq 
The above identity can be used to eliminate all prefactors 
involving $u_n$ wherever 
they appear, so that the interesting dynamical relations that remain will only involve 
coefficients of $P_n$ and $Q_n/u_n$. Note also that, from (\ref{parq}), the recursion 
breaks down if $d_n$ vanishes at some stage. 

There is a natural geometrical interpretation of the iteration that produces the 
continued fraction expansion (\ref{cfrac}) for (\ref{pick}), which goes back to  
work of Adams and Razar on the elliptic case \cite{ar}, and 
was further generalized by Bombieri and Cohen in the setting of 
Pad\'e approximation of functions on algebraic curves of general type \cite{bombcoh}. 
Let us consider the function 
\beq\label{fdef} 
G= Y_0-a_0 = Y_1^{-1}, 
\eeq 
and denote by $x_0^{(1)},\ldots,x_0^{(\rg)}$ 
and $x_1^{(1)},\ldots,x_1^{(\rg)}$ the roots of 
$Q_0(X)$ and $Q_1(X)$, respectively. If we also 
set 
$$
y^{(j)}_n =P_n(x^{(j)}_n), \qquad j=1,\ldots,\rg 
$$ 
for $n=0,1$, then under the Abel map each of the 
degree $\rg$ divisors 
$$ 
D_0 = (x_0^{(1)},y_0^{(1)})+ \cdots + (x_0^{(\rg)},y_0^{(\rg)}), 
\qquad 
D_1 = (x_1^{(1)},y_1^{(1)})+ \cdots + (x_1^{(\rg)},y_1^{(\rg)})
$$ 
corresponds to a point on the Jacobian variety 
$\mathrm{Jac}({\cal C})\cong\mathrm{Sym}^\rg({\cal C})$, 
identified with the $\rg$-fold symmetric product of the curve \cite{mumford}. 
From its expression as $Y_0-a_0$, the poles of $G$ lie at the points 
$(x_0^{(1)},y_0^{(1)}), \ldots, (x_0^{(\rg)},y_0^{(\rg)})$ 
and $\infty_2$, and it vanishes precisely at 
$(x_1^{(1)},y_1^{(1)}), \ldots, (x_1^{(\rg)},y_1^{(\rg)})$ 
and $\infty_1$, where $Y_1$ has poles. Therefore $D_0$ and 
$D_1$ are related by the linear equivalence 
\beq\label{lineq} 
D_1\sim_{\mathrm{lin}} D_0 +\infty_2-\infty_1, 
\eeq 
and by the same argument the shift $n\to n+1$ in each line of the 
continued fraction is equivalent to a translation in the Jacobian 
by the divisor class of $\infty_2-\infty_1$.

\section{Lax pair and nonlinear system} \label{laxnonl} 
\setcounter{equation}{0}

The recursion (\ref{rec}) and the relations (\ref{eig}) can be reformulated in 
terms of a linear system, which makes their structure much more transparent. 
To do this, it suffices to introduce projective 
coordinates, setting  
$$Y_n=\psi_n/\phi_n, \qquad  
\mathbf{\Psi}_n=(\psi_n,  \phi_n)^T, $$
and substituting into (\ref{eig}), which  leads to the eigenvalue problem 
\beq\label{eigval} 
{\bf L}_n(X)\, \mathbf{\Psi}_n=Y \mathbf{\Psi}_n 
\eeq 
for the Lax matrix 
\beq\label{lax} 
{\bf L}_n= \left(\begin{array}{cc} P_n & Q_{n-1} \\ 
Q_n & -P_n \end{array}\right). 
\eeq 
Upon substituting the ratio of the projective coordinates into (\ref{rec}), 
and fixing an arbitrary multiplier, the fractional linear relation between 
$Y_n$ and $Y_{n+1}$ separates into two linear equations, which can be 
written as 
\beq\label{mshift} 
{\bf M}_n(X)\mathbf{\Psi}_{n+1}=  
\mathbf{\Psi}_n, 
\eeq 
where, taking the standard formula from the classical theory of continued fractions 
for real numbers, we may set 
\beq\label{em} 
{\bf M}_n =\left(\begin{array}{cc} a_n & {1} \\ 
1 & 0 \end{array}\right) .
\eeq
(The choice of multiplier means that the matrix  ${\bf M}_n$ is only defined up to 
overall scaling ${\bf M}_n\to \la_n {\bf M}_n$  for  some arbitrary $n$-dependent quantity $\la_n$.) 

The compatibility of the linear system consisting of (\ref{eigval}) and 
(\ref{mshift}) is the discrete Lax equation  
\beq\label{dlax} 
{\bf L}_{n}{\bf M}_n
= {\bf M}_n{\bf L}_{n+1}, 
\eeq 
which produces two non-trivial conditions, namely 
\beq \label{compat} 
\begin{array}{rcl}
(X+v_n)(P_{n+1}-P_n)& = & \frac{1}{2}u_n(Q_{n-1}-Q_{n+1})\\
(X+v_n)Q_n & = & \frac{1}{2}u_n(P_{n+1}+P_n), 
 \end{array} 
\eeq 
where we have substituted the expression (\ref{parq}) for the partial quotient $a_n$. 
Note also that, because it is equivalent to conjugation by the nonsingular matrix ${\bf M}_n$, 
(\ref{dlax}) is an isospectral evolution, preserving the spectral curve 
$\det (Y\mathbf{1}-{\bf L}_n(X))=0$, which reproduces the formula 
(\ref{dety}) for all $n$.   

Now, from the terms in the continued fraction, we can expand 
\beq\label{pqexp} 
P_n(X)=A(X) +\sum_{j=0}^{\rg-1}\pi^{(j)}_nX^j, 
\qquad 
Q_n(X) = u_n \, \Big( X^\rg +\sum_{j=0}^{\rg-1}\rho^{(j)}_nX^j\Big), 
\eeq 
with two particular coefficients being specified as 
\beq \label{topPQ}
\pi^{(g-1)}_n = 2d_n, \qquad \rho^{(g-1)}_n =-v_n
\eeq
in terms of the notation used previously. Upon substituting (\ref{pqexp}) into 
(\ref{compat}), the terms involving $u_n$ can be replaced using (\ref{ud}) 
in the first relation, and cancelled from the second relation, to 
yield a set of recurrences for the coefficients $\pi^{(j)}_n$, $\rho^{(j)}_n$ 
in (\ref{pqexp}), namely 
\beq\label{recs1} 
(X+v_n)\sum_{j=0}^{\rg-1}(\pi^{(j)}_{n+1}-\pi^{(j)}_n)X^j  
=2(d_{n+1}-d_n)X^\rg + 2\sum_{j=0}^{\rg-1}(d_{n+1}\rho^{(j)}_{n+1}-d_n\rho^{(j)}_{n-1})X^j,   
\eeq 
\beq\label{recs2} 
(X+v_n)\Big(X^\rg+ \sum_{j=0}^{\rg-1}\rho^{(j)}_{n}X^j  \Big) 
=A(X)+\frac{1}{2}\sum_{j=0}^{\rg-1}(\pi^{(j)}_{n+1}+\pi^{(j)}_{n})X^j.   
\eeq  
Let us introduce the $\rg$-tuples 
$$ 
\bm{\pi}_n =( \pi^{(0)}_{n},\ldots,\pi^{(\rg-1)}_{n}), 
\qquad 
\bm{\rho}_n =(\rho^{(0)}_{n},\ldots,\rho^{(\rg-1)}_{n}) 
$$ 
of affine coordinates.
Due to (\ref{aform}) and (\ref{topPQ}), the coefficients at order $X^\rg$ 
in (\ref{recs1}) and at orders $X^{\rg+1}$ and  $X^\rg$ in (\ref{recs2}) 
provide only tautologies, so that altogether there are $2\rg$ non-trivial relations 
between the components of $\bm{\pi}_{n},\bm{\pi}_{n+1}$,
and $\bm{\rho}_{n},\bm{\rho}_{n\pm1}$.
To be precise, these relations mean that the $2\rg$ quantities 
$\bm{\pi}_{n+1}$ and $\bm{\rho}_{n+1}$ 
can be calculated as rational functions of the 
components of $\bm{\pi}_{n}$,
$\bm{\rho}_{n-1}$, and  $\bm{\rho}_{n}$, 
which (together with the relation (\ref{ud}) for the prefactors) shows how 
the entries of ${\bf L}_{n+1}$ are determined from those of ${\bf L}_n$. 
Similarly, in the reverse direction $n+1\to n$, these relations mean that the 
entries of  ${\bf L}_n$ can be obtained as rational functions of the entries 
of ${\bf L}_{n+1}$. 

In the above form, the map corresponding to the shift $n\to n+1$ from one line 
of the continued fraction to the next can be interpreted as a discrete dynamical system, 
where 
(ignoring the prefactors $u_n$) 
this can be viewed as a birational map  
$
(\bm{\pi}_n,\bm{\rho}_{n-1},\bm{\rho}_n) \mapsto 
(\bm{\pi}_{n+1},\bm{\rho}_{n},\bm{\rho}_{n+1})
$
in dimension $3\rg$. However, at the 
expense of introducing more parameters, one can use the equation for the spectral 
curve (\ref{dety}) to eliminate $\rg$ coordinates and rewrite this in terms 
of a birational map in dimension $2\rg$. In particular, 
in  the explicit formula (\ref{specrel}) the leading order ($X^{2\rg}$) term gives 
(\ref{ud}), while the coefficients at each order  from
$X^{2\rg-1}$ down to $X^\rg$   can be used to rewrite $\rho^{(0)}_{n-1},\ldots,\rho^{(g-1)}_{n-1}$ 
in terms of the components of $\bm{\pi}_n$, $\bm{\rho}_n$, 
as well as the coefficients in appearing in $A(X)$, 
and also $u$, the leading coefficient of $R(X)$ on the left-hand side. The remaining $\rg$ coefficients 
of $R(X)$, which appear at orders $X^j$, $j=0,\ldots,\rg-1$, can then be written as rational functions 
of the components of $\bm{\pi}_n$, $\bm{\rho}_n$ and the other parameters, 
and these $\rg$ quantities are independent of $n$. In this way, we arrive 
at a birational map 
\beq \label{phi} 
\varphi: \qquad (\bm{\pi}_n,\bm{\rho}_n) \mapsto 
(\bm{\pi}_{n+1},\bm{\rho}_{n+1})
\eeq 
which has $\rg$ conserved quantities. In fact there is more that one can say about this map: it turns 
out to be symplectic, and integrable in the sense of a suitable discrete analogue of Liouville's theorem 
\cite{bruschi,maeda,veselov}. 

\begin{thm}\label{int} 
The birational map (\ref{phi}) corresponding to the iteration for the   
continued fraction expansion (\ref{cfrac}) of the function (\ref{pick}) 
is an integrable symplectic map in dimension $2\rg$. 
\end{thm} 

Observe that the expression (\ref{dety}) is symmetrical in $Q_{n-1}$ and $Q_n$, so one can just as well 
use it eliminate the components of $\bm{\rho}_n$ from (\ref{recs1}) and (\ref{recs2}), to obtain 
a birational map 
\beq \label{phihat} 
\hat{\varphi}: \qquad (\bm{\pi}_{n},\bm{\rho}_{n-1}) \mapsto 
(\bm{\pi}_{n+1},\bm{\rho}_{n}).
\eeq
Clearly the latter map is conjugate to $\varphi$, in the sense that there is a birational transformation $\chi$ 
such that $\hat{\varphi}= \chi^{-1}\circ\varphi\circ\chi$, and the above 
theorem applies equally well to $\hat\varphi$. 
The  general proof this theorem is given in section \ref{poisson}, where we make use of a 
Poisson structure for the Lax matrices (\ref{lax}). For now, we just give explicit details 
for $\rg=1$ and 2.  

\begin{exa}\label{elliptic} 
{\bf The case $\rg=1$: } In the genus one case, following \cite{vdp}, we write  
\beq\label{g1co} 
A(X) = X^2+f, \qquad P_n=A(X)+2d_n, \qquad Q_n =u_n(X-v_n), \qquad R(X)=u(X-v)
\eeq 
for arbitrary parameters $f,u,v$ defining the quartic curve 
$$Y^2=(X^2+f)^2+4u(X-v)
$$ in the $(X,Y)$ plane.
There are only two non-trivial relations from (\ref{recs1}) and   (\ref{recs2}), 
given by  
\beq\label{g1rels}
d_{n+1}+d_n +v_n^2+f=0, \qquad 
d_{n+1}(v_{n+1}+v_n)=d_n(v_n+v_{n-1}), 
\eeq 
which define a birational map in 3 dimensions, that is 
$$ 
(d_n, v_{n-1},v_n) \mapsto (d_{n+1},v_n,v_{n+1}) =\left(-d_n -v_n^2-f,v_n, 
-v_n -\frac{d_n(v_n+v_{n-1})}{(d_n + v_n^2+f)}\right) .
$$
However, using the equation for the curve, and removing an overall factor of 4, 
the formula (\ref{specrel}) becomes 
\beq\label{g1}
u(X-v) =d_n(X^2+f)-d_n(X-v_{n-1})(X-v_n)+d_n^2. 
\eeq 
The first non-trivial relation, at order $X$, gives 
$$ 
v_{n-1}=-v_n +\frac{u}{d_n}, 
$$ 
which allows  $v_{n+1}$ to be rewritten as a function of $d_n$, $v_n$ and the parameters $f,u$. 
Hence, making use of (\ref{g1rels}),  this yields a map in the plane, that is 
\beq\label{g1map} 
\varphi: \qquad 
(d_n,v_n)\mapsto(d_{n+1},v_{n+1}) = \left(-d_n -v_n^2-f, 
-v_n -\frac{u}{(d_n + v_n^2+f)}\right) .
\eeq 
The above map preserves the symplectic form 
$$ 
\om = \rd d_n \wedge \rd v_n, 
$$ 
or in other words 
$\varphi^* \om = \rd d_{n+1}\wedge \rd v_{n+1} = \om$. 
Furthermore,  the lowest order ($X^0$) term in (\ref{g1}) provides the relation 
$$ 
-uv=d_n(f+d_n-v_{n-1}v_n), 
$$ 
so replacing $v_{n-1}$ as before and setting $H=-uv$ we see that 
\beq\label{g1ham}
H= d_nv_n^2-uv_n +d_n^2+fd_n
\eeq
is a conserved quantity for $\varphi$, so it is an integrable symplectic map in two dimensions. 
\end{exa} 

\begin{figure} \centering 
\includegraphics[width=11cm,height=11cm,keepaspectratio]{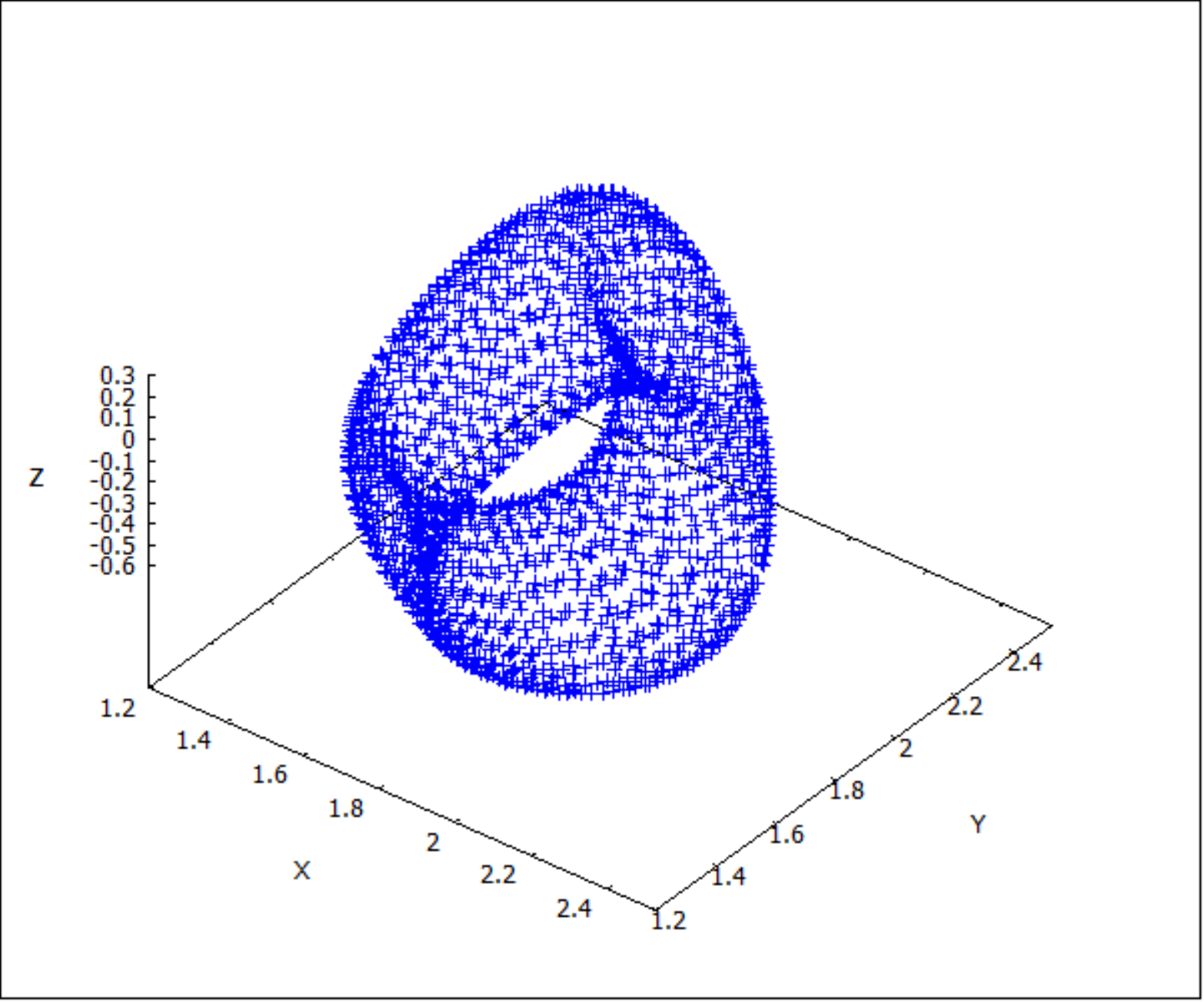}
\caption{\small{The 3D projection $(d_{n-1},d_n,v_{n-1})$ of 2000 points on 
the orbit of the 4D map   
  (\ref{g2pair}) with $f=-5$, $g=u=-1$ and 
$(d_0,d_1,v_0,v_1)=(5/4,2,-1/2,0)$.  }}
\label{f1}
\end{figure}

\begin{exa}\label{genus2}
{\bf The case $\rg=2$: } In the genus two case,  adopting the  notation in \cite{vdp2}, we 
write \beq\label{g2co} 
\begin{array}{l} 
A(X) = X^3+fX+g, \quad P_n=A(X)+2d_n(X+e_n), \\ 
Q_n =u_n(X^2-v_nX+w_n), \quad R(X)=u(X^2-vX+w)
\end{array} \eeq 
for arbitrary parameters $f,g,u,v,w$ defining the sextic curve 
\beq\label{csextic} {\cal C}: \qquad
Y^2=(X^3+fX+g)^2+4u(X^2-vX+w). 
\eeq 
From  (\ref{recs1}) and   (\ref{recs2}) there are four relations that 
define a birational map in 6 dimensions, namely 
\beq\label{rel1}
d_{n+1}(e_{n+1}+v_{n+1}+v_n)= d_n(e_n+v_n+v_{n-1}),  
\eeq 
\beq\label{rel2} 
v_n(d_{n+1}e_{n+1}-d_ne_n)=d_{n+1}w_{n+1}-d_nw_{n-1}, 
\eeq 
\beq\label{rel3} 
d_{n+1}+d_n+f=w_n-v_n^2, 
\eeq
\beq\label{rel4} 
d_{n+1}e_{n+1}+d_ne_n+g=v_nw_n. 
\eeq
To be precise, we have the map 
$(d_n,e_n,v_{n-1},v_n,w_{n-1},w_n)\mapsto (d_{n+1},e_{n+1},v_{n},v_{n+1},w_{n},w_{n+1})$, 
where (\ref{rel3}) is used to obtain $d_{n+1}$, and then (\ref{rel4}) produces an expression for 
$e_{n+1}$, which allows $v_{n+1}$ and $w_{n+1}$ to be calculated from (\ref{rel1}) and (\ref{rel2}), 
respectively.  
In order to obtain a map in 4 dimensions, one can use 
(\ref{specrel}) to eliminate $v_{n-1}$ and $w_{n-1}$, giving (\ref{phi}), or   
instead eliminate $v_{n}$ and $w_{n}$, to obtain (\ref{phihat}); 
here we take the latter option. To be precise, compared with the 
quantities used in (\ref{phihat}) we have made an invertible   change 
of coordinates, that is  
$(\pi_n^{(0)}, \pi_n^{(1)},\rho_{n-1}^{(0)}, \rho_{n-1}^{(1)}) 
=(d_ne_n, d_n, w_{n-1},-v_{n-1})$, but by a slight abuse of  
of notation we will use the same symbol $\hat\varphi$ to denote 
the map that describes the shift $n\to n+1$ in terms of 
$(d_n,e_n,v_{n-1},w_{n-1})$.   There are four non-trivial relations 
coming from (\ref{specrel}), given by 
\beq\label{nrel1} 
d_n(e_n +v_n+v_{n-1})=0, 
\eeq 
\beq\label{nrel2} 
u=d_n(d_n-v_nv_{n-1}-w_n-w_{n-1}+f), 
\eeq 
\beq\label{nrel3} 
-uv=d_n(2d_ne_n+v_nw_{n-1}+v_{n-1}w_n +fe_n+g), 
\eeq 
\beq\label{nrel4} 
uw=d_n(d_ne_n^2-w_nw_{n-1}+ge_n). 
\eeq  
Using (\ref{nrel1}) and (\ref{nrel2}), together with (\ref{rel3}) and (\ref{rel4}), we find the map 
\beq\label{g2phihat} 
\hat{\varphi}: \qquad 
(d_n,e_n,v_{n-1},w_{n-1}) \mapsto 
 (d_{n+1},e_{n+1},v_{n},w_{n}) ,  
\eeq 
where the shifted variables are given explicitly in terms of the previous ones by 
$$ \begin{array}{rcl} 
d_{n+1}& = & -e_n^2-e_nv_{n-1}-w_{n-1}-\frac{u}{d_n}=:-d_n^{-1}{\cal D} , \\  
e_{n+1}& =& {\cal D}^{-1}\left[v_{n-1}\Big(d_n^2+d_n(e_n+v_{n-1})^2-d_nw_{n-1}+fd_n-u\Big) 
+e_n\Big(2d_n^2-d_nw_{n-1}+fd_n-u\Big)+gd_n\right] , 
\\
v_n& =&  
-v_{n-1}-e_n,\\
w_n & =& -w_{n-1}+v_{n-1}^2+v_{n-1}e_n+d_n 
-ud_n^{-1}+f. 
\end{array}
$$ 
This four-dimensional map preserves a nondegenerate Poisson 
bracket, given by 
$$ 
\{\,d_n,e_n\,\} = \{\,d_n,v_{n-1}\,\}=\{\,v_{n-1},w_{n-1}\,\}=0, 
$$ 
$$ 
\{\,d_n,w_{n-1}\,\}=-1, \quad \{\,e_n,v_{n-1}\,\}=\frac{1}{d_n}, 
\quad
\{\,e_n,w_{n-1}\,\}=\frac{v_{n-1}+e_n}{d_n}, 
$$ 
hence it is symplectic. (As we shall see in section \ref{poisson}, 
in a different set of coordinates, with $e_n$  replaced by 
$\pi_n^{(0)}=d_ne_n$, this bracket has only linear and constant terms.)  
Upon eliminating $v_{n-1}$ and $w_{n-1}$, one can also rewrite this as 
a map $(d_{n-1},d_n,v_{n-1},v_n)\mapsto (d_n,d_{n+1},v_{n},v_{n+1})$, 
so that it takes a simpler form as a pair of coupled recurrence relations of 
second order, that is 
\beq\label{g2pair} 
\begin{array}{rcl}
d_{n+1}+d_{n}+d_{n-1} +u/d_n 
+v_n^2+v_nv_{n-1}+v_{n-1}^2+f & =& 0 , \\
(2v_n+v_{n-1})d_n +
(2v_n+v_{n+1})d_{n+1} +v_n^3+fv_n-g &=&0, 
\end{array} 
\eeq 
and in these coordinates (up to an overall choice of scaling) the symplectic form is  
$$
\omega =\mathrm{d}d_{n-1} \wedge  \mathrm{d}d_n 
+ (2v_{n-1}+v_n)\, 
\mathrm{d}v_{n-1} \wedge  \mathrm{d}d_n
+
d_n\,\mathrm{d}v_{n-1} \wedge  \mathrm{d}v_n 
 . 
$$ 
By construction, both of the quantities $H_1=-uv$, $H_2=uw$ defined 
by (\ref{nrel3}) and (\ref{nrel4}) are conserved, 
and it  can  be verified directly  
that $\{\,H_1,H_2\,\}=0$, hence the map is integrable in the Liouville sense. 
A particular orbit of this map is plotted in Fig.~\ref{f1}.
\end{exa} 

\section{
Orthogonal polynomials and Hankel determinants}  \label{hankels} 
\setcounter{equation}{0}  

After removing the first term $a_0$, the continued fraction expansion  (\ref{cfrac}) 
becomes 
\beq\label{gfrac} 
G=
\cfrac{ 1}{ a_1 +\cfrac {1}{ a_2 +\cdots} }  
=\cfrac{ 1}{ 2u_1^{-1}(X+v_1) +\cfrac {1}{ 2u_2^{-1}(X+v_2) +\cdots} }. 
\eeq A continued fraction of this form, where each 
partial quotient $a_j = a_j(X)$ is a linear function of $X$, is called a
J-fraction \cite{wall}. 
If we multiply the main numerator and denominator by $u_1/2$ and apply (\ref{ud}), 
then this can be rewritten in its more classical 
form, that is 
\beq\label{gcla} 
G
=\cfrac{ -2u_0^{-1}d_1}{X+v_1 -\cfrac {d_2}{ X+v_2-\cfrac{d_3}{X+v_3 -\cdots}} },  
\eeq
from which we see that, up to an overall prefactor of $-2/u_0$, it is completely determined 
by the quantities $d_n$ and $v_n$. 
In this section we apply standard results 
on J-fractions and associated orthogonal polynomials, which 
lead directly to formulae for the quantities $d_n$ and $v_n$ in terms 
of ratios of Casorati determinants, and Hankel determinants in 
particular.   
This generalizes certain results obtained for $\rg=1$ by Chang, Hu and  Xin \cite{chx} 
to hyperelliptic curves of any genus $\rg$.   

In the neighbourhood of the point $\infty_1\in{\cal C}$, the function $G$ defined by 
(\ref{fdef}) 
has the series expansion 
\beq\label{gseries} 
G=\sum_{j\geq 1}s_{j-1} \, X^{-j} = s_0 X^{-1}+s_1X^{-2}+s_2X^{-3}+\ldots,  
\eeq 
which can be used to define a linear functional $<\cdot >: {\cal F}\to \C$ 
according to 
\beq\label{func} 
<\Phi>=\frac{1}{2\pi\ri}\oint_{\infty_1}\Phi \,G\,\rd X, 
\eeq 
for any function $\Phi\in{\cal F}$, where the integral is taken along any sufficiently small closed contour 
around $\infty_1$ (anticlockwise in the $X$ plane) that does not encircle the poles of $G$ at 
$(x_0^{(1)},y_0^{(1)}), \ldots, (x_0^{(\rg)},y_0^{(\rg)})$ 
and $\infty_2$. In other words, $G$ can be regarded as 
a moment generating function with moments 
\beq\label{moms} 
<X^j> = s_{j}, \qquad j=0,1,2,\ldots, 
\eeq 
although in general, for complex $s_j$,  this may not be associated  
with some positive measure. 
The linear functional (\ref{func}) 
also defines a scalar product between any pair 
of functions $\Phi$, $\Psi$, that is 
\beq\label{prod}
<\,\Phi,\Psi\,> := <\Phi\Psi>.
\eeq  

The convergents of (\ref{gfrac}) are  the sequence of rational functions of $X$ given by 
$$ 
\frac{p_0}{q_0}=0, \quad \frac{p_1}{q_1}=\frac{1}{a_1},\quad  
   \frac{p_2}{q_2}=
\cfrac{ 1}{ a_1 +\cfrac {1}{ a_2 } } , 
\quad \frac{p_3}{q_3}=
\cfrac{ 1}{ a_1 +\cfrac {1}{ a_2 +\cfrac{1}{a_3}} } , \quad \ldots , 
$$ 
where by convention one can take  
\beq \label{inits1} 
p_0=0, \qquad q_0=1,\eeq 
followed by  
\beq \label{inits2}
p_1=s_0 = \frac{u_1}{2}, \qquad 
q_1 =X+v_1 = X-\frac{s_1}{s_0},  
\eeq 
and for all $j\geq 1$ 
$p_j(X)$ and $q_j(X)$ are polynomials of degrees $j-1$ and $j$ in $X$, respectively, where  without loss 
of generality each $q_j$ is taken to be monic. The recursion for the convergents 
is essentially controlled by the same matrix (\ref{em}) that appears in the 
classical theory of continued fraction expansions of numbers  in $\R$. 
However, 
the entries of this matrix must be scaled in order to ensure that 
the $q_j$ are monic, 
making use of (\ref{ud}) to remove dependence on the prefactors $u_j$, 
to yield  
$$ 
\left(\begin{array}{cc} p_n & q_n  \\ 
p_{n-1} & q_{n-1} \end{array}\right) =\left(\begin{array}{cc} X+v_n & -d_n  \\ 
1 & 0 \end{array}\right) \left(\begin{array}{cc} p_{n-1} & q_{n-1} \\ 
p_{n-2} & q_{n-2} \end{array}\right) 
,
$$ 
which is equivalent to the three-term linear recurrence relation 
\beq\label{3term} 
q_n=(X+v_n)\, q_{n-1}-d_n\, q_{n-2}, \qquad n\geq 2, 
\eeq 
for the sequence of polynomials $(q_n)$, and the same for 
the sequence $(p_n)$, with the initial conditions (\ref{inits1}) 
and (\ref{inits2}). From the linear recurrence it is clear that $p_n/s_0$ is monic for each $n$.     

Due to the fact that each $a_j=\lc Y_j \rc$   is linear  in $X$, it is 
straightforward to see by induction that the $n$th 
approximant $p_n/q_n$  satisfies  
$$ 
G-\frac{p_n}{q_n} = \cfrac{ 1}{ a_1 +\cfrac {1}{ \cdots +\cfrac{1}{a_n+Y_{n+1}^{-1}}} } 
- \cfrac{ 1}{ a_1 +\cfrac {1}{ \cdots +\cfrac{1}{a_n}} } 
= O\left(X^{-(2n+1)}\right), 
$$ 
and since $q_n$ has degree $n$ this implies that 
\beq\label{linap} 
p_n - Gq_n = O(X^{-(n+1)}), \qquad n=0,1,2,\ldots . 
\eeq 
Given the requirement on the degrees of $p_n$ and $q_n$, 
by considering the terms at orders $X^{n-1},X^{n-2},\ldots,X^{-n}$, 
the equation (\ref{linap}) 
provides $n$ linear equations that determine the non-trivial coefficients of 
$q_n$   in terms of the coefficients $s_j$ in (\ref{gseries}), 
and a further $n$ linear equations for the non-trivial 
coefficients of $p_n$ in terms of those of $q_n$ and the $s_j$. 
This leads to a standard formula for $q_n$, 
given explicitly in determinantal form 
as  
\beq\label{qform}
q_n(X) = \frac{1}{\Delta_n} \, 
\left| \begin{array}{ccccc} 
s_0     & s_1 & \cdots & s_{n-1}  & s_{n} \\ 
s_1     &         & \iddots &       &  \vdots       \\
\vdots & \iddots &          &       &     \vdots    \\        
s_{n-1}     &  \cdots & \cdots        & \cdots & s_{2n-1} \\ 
1        &  X         & X^2 & \cdots & X^n 
\end{array} 
\right| , 
\eeq 
where $\Delta_n$ is the $n\times n$ Hankel determinant
\beq\label{hankel} 
\Delta_n = 
\left| \begin{array}{cccc} 
s_0     & s_1 & \cdots & s_{n-1}   \\ 
s_1     &         & \iddots &        \vdots       \\
\vdots & \iddots &          &           \vdots    \\        
s_{n-1}     &  \cdots & \cdots & s_{2n-2}  
\end{array} 
\right| =\det (s_{i+j-2})_{i,j=1,\ldots, n}. 
\eeq  
For what follows, we also need to introduce another determinant of Casorati type, 
obtained from the coefficients $s_j$ in (\ref{gseries}) by shifting the last column of 
the Hankel matrix, namely 
\beq\label{cas}
 \Delta_n^* = 
\left| \begin{array}{ccccc} 
s_0     & s_1 & \cdots & s_{n-2} & s_{n}    \\ 
s_1     &         & \iddots &    \vdots        & \vdots       \\
\vdots & \iddots &          &  \vdots       &   \vdots    \\ 
s_{n-2}     &  \cdots & \cdots & s_{2n-4} &  s_{2n-2}   \\
s_{n-1}     &  \cdots & \cdots & s_{2n-3} &  s_{2n-1}  
\end{array} 
\right|.
\eeq 
We refer to the latter as a shifted Hankel determinant. 
By convention we set $\Delta_0=1$ and $\Delta_0^*=0$.

We can now state our main result about Hankel determinants and orthogonal 
polynomials.
 
\begin{thm}\label{hankt} 
The quantities $d_n$, $v_n$ that appear under iteration of the J-fraction expansion 
(\ref{gcla}) of the 
function (\ref{fdef}), which provide the two components (\ref{topPQ}) 
of the iterates of the birational map (\ref{phi}), are given in 
terms of  Hankel and shifted Hankel determinants by 
\beq\label{hankform} 
d_n = \frac{\Delta_n\Delta_{n-2}}{\Delta_{n-1}^2} \quad\mathrm{for}\,\,n\geq 2, 
\qquad v_n = \frac{\Delta_{n-1}^*}{ \Delta_{n-1}}
- \frac{\Delta_{n}^*}{ \Delta_{n}} \quad \mathrm{for}\,\,n\geq 1, 
\eeq 
where the entries 
are determined recursively by $s_0=u_1/2$ and,  for $j\geq 1$, 
\beq\label{sj} 
s_j = -\sum_{i=2}^{\rg+1}(k^{(\rg+1-i)} + \pi_1^{(\rg+1-i)}) s_{j-i} 
+\frac{d_1}{s_0}\left(\sum_{i=0}^{j-2} s_is_{j-2-i} 
+\sum_{\ell=3}^{\rg+2}\sum_{i=0}^{j-\ell}\rho_0^{(\rg+2-\ell)} s_i s_{j-\ell-i}\right) 
+\bar{s}_j
\eeq 
where 
$$ 
\bar{s}_j = 
\begin{cases}   
s_0 \rho_1^{(\rg-j)} & \mathrm{for }\,\, 1\leq j\leq \rg, \\ 
0 &\mathrm{for } \,\, j\geq \rg+1. 
\end{cases} 
$$ 
Furthermore, the polynomials $q_j$ that appear as the denominators 
of the convergents of the J-fraction are orthogonal with respect 
to the scalar product (\ref{prod}) associated with the series (\ref{gseries}), 
that is 
\beq\label{orthog}
<\,q_n,q_m\,> = h_n\, \delta_{nm},  
\eeq 
where 
$$ 
h_n =\frac{\Delta_{n+1}}{\Delta_n}. 
$$ 
\end{thm} 
\begin{prf} 
The formulae in (\ref{hankform}) are classical expressions 
for the coefficients appearing in the linear recurrence  (\ref{3term}) 
for orthogonal polynomials. A direct proof is obtained by substituting 
(\ref{qform}) into the three-term recurrence and expanding in powers of $X$: 
the determinantal expression for $v_n$ appears immediately at order $X^n$, 
while at order $X^{n-1}$ one finds a formula for $d_n$ as a combination of four 
terms that can be condensed into a single ratio 
by applying various identities for determinants (which are collected in 
an appendix for completeness). To find the entries of the Hankel determinants 
recursively, note from (\ref{y1}) that $G=1/Y_1=(Y-P_1)/Q_0$, hence 
$Y=Q_0G+P_1$ and 
$$ 
Y^2=Q_0^2G^2+2Q_0P_1G+P_1^2=P_1^2+Q_0Q_1, 
$$ 
using (\ref{dety}), so that $G$ satisfies the quadratic equation 
\beq\label{genfn} 
P_1G+\frac{1}{2}Q_0 G^2=\frac{1}{2}Q_1. 
\eeq 
Upon substituting the series (\ref{gseries}) for $G$ and expanding 
the polynomial coefficients in (\ref{genfn})  with the notation  in 
(\ref{aform}) and (\ref{pqexp}), and making use of (\ref{ud}) to replace $u_0$, 
the recursion (\ref{sj}) results. 
Finally, the orthogonality of the sequence of polynomials 
$(q_n)$ follows by a standard inductive argument using 
the three-term recurrence, also making use of the 
moments (\ref{moms}) with (\ref{qform}) to expand 
$\Delta_n^2<q_n^2>$ as a sum of $n+1$ products of 
determinants, only one of which is non-vanishing, which 
yields $<\,q_n,q_n\,>=\Delta_{n+1}/\Delta_n=h_n$.    
\end{prf} 

\begin{exa}\label{elliphankel} 
{\bf The recursion for moments in the elliptic case:} 
For $\rg=1$ we use the same notation as in Example \ref{elliptic}, 
and for the recursion (\ref{sj}) we have 
$$s_0=u_1/2, \qquad  
s_1=-s_0v_1$$  and  
\beq\label{elliprec} 
s_j = -(2d_1+f)s_{j-2} +s_0^{-1}d_1\left( 
\sum_{i=0}^{j-2} s_i s_{j-2-i} -v_0 \sum_{i=0}^{j-3} s_i s_{j-3-i}\right) 
\eeq 
for all $j\geq 2$. 
To illustrate this with a particular numerical example, let us pick the curve 
\beq\label{s4quart} 
Y^2=(X^2-3)^2 -4(X+2), 
\eeq  
so that $f=-3$, $u=-1$, $v=-2$, and 
set $u_0=-2$,  $d_0=1$, $v_0=-1$. This fixes the function 
\beq\label{gs4} 
G=Y_0-a_0=
\frac{1+\frac{1}{2}(X^2-3-Y)}{X+1}
= {X}^{-1}+
{2}X^{-3}+
X^{-4}+
6X^{-5}+\cdots .
\eeq
Then 
 (\ref{g1map}) produces the values $d_1=1$, $v_1=0$, which gives 
$s_0=u_1/2=-2d_1/u_0=1$ and $s_1=-s_0v_1=0$, and  
the coefficients of the series expansion  
(\ref{gs4})  
are obtained from the particular recurrence 
$$ s_j = s_{j-2} +
\sum_{i=0}^{j-2} s_i s_{j-2-i} + \sum_{i=0}^{j-3} s_i s_{j-3-i}, \quad j\geq 2 .
$$
This produces the sequence of moments 
$$ 
1,0,2,1,6,7,24,41,115,236,613,1380,\ldots, 
$$ and 
the corresponding sequence of Hankel determinants 
begins with 
\beq\label{haseq} 
\Delta_0=\Delta_1=1,  
\quad \Delta_2
=\left|\begin{array}{cc} 1 & 0 \\ 0 & 2 \end{array}\right|=2, 
\quad \Delta_3= \left|\begin{array}{ccc} 1 & 0 & 2 \\ 0 & 2 & 1 \\ 2 & 1 & 6 \end{array}\right|=3, 
\quad \Delta_4= \left|\begin{array}{cccc} 1 & 0 & 2 & 1 \\ 0 & 2 & 1 & 6 \\ 2 & 1 & 6 & 7 \\ 1& 6 & 7 & 24 \end{array}\right|=7, \ldots, 
\eeq 
which should remind the reader of (\ref{s4seq}). 
\end{exa} 

\begin{exa}\label{g2hankel} 
{\bf The recursion for moments in genus 2:}
With the notation of Example \ref{genus2}, 
the recursion (\ref{sj}) for $\rg=2$ has initial values  
$$s_0=u_1/2, \qquad  
s_1=-s_0v_1, \qquad s_2 = s_0(w_1-d_1-f),$$  and 
subsequent coefficients in the series expansion (\ref{gseries}) are 
determined for   $j\geq 3$ by 
\beq\label{g2rec} 
s_j = -(2d_1+f)s_{j-2} -(2d_1e_1+g)s_{j-3} 
+s_0^{-1}d_1\left( 
\sum_{i=0}^{j-2} s_i s_{j-2-i} -v_0 \sum_{i=0}^{j-3} s_i s_{j-3-i}
+w_0 \sum_{i=0}^{j-4} s_i s_{j-4-i}\right) 
. 
\eeq 
As a particular example, consider the curve 
$$ 
Y^2=(X^3-5X-1)^2-4(X^2+2X+3), 
$$ 
with $f=-5$, $g=u=-1$, $v=-2$, $w=3$, and choose 
$u_0=-4$, $d_0=5/4$, $e_0=3/5$, $v_0=-1/2$, $w_0=-3/2$, 
corresponding to the orbit of the map (\ref{g2pair}) plotted in 
Fig.~\ref{f1}, 
which gives 
$$ 
G=Y_0-a_0 = \frac{X+\frac{1}{2}+\frac{1}{4}(X^3-5X-1-Y)}{X^2+\frac{1}{2}X-\frac{3}{2}} 
=X^{-1}+2X^{-3}+7X^{-5}+2X^{-6}+\cdots, 
$$ 
and the recursion for the coefficients (moments)  in the above 
expansion   is 
$$ 
s_j=s_{j-2}-s_{j-3} 
+ 2 \sum_{i=0}^{j-2} s_i s_{j-2-i} + \sum_{i=0}^{j-3} s_i s_{j-3-i}
-3 \sum_{i=0}^{j-4} s_i s_{j-4-i}, \quad j\geq 3
$$  
with initial values $s_0=1$, $s_1=0$, $s_2=2$. In this 
case,  
the sequence of moments $(s_j)$ begins with  
$$ 
1,0,2,0,7,2,31,21,159,168,900,1246,5455,9040,34731,65328,\ldots , 
$$ 
yielding the corresponding sequence of Hankel determinants 
\beq\label{g2hank} 
(\Delta_j): \quad 
1,1,2,6,31,319,5810, 147719, 8526736, \ldots ,  
\eeq 
as well as the sequence of shifted Hankel determinants 
$$ 
(\Delta_j^*): \quad 
0,0,0,4,16,200,6987,161401,11022617,\ldots . 
$$ 
\end{exa} 

The map (\ref{phi}), or equivalently (\ref{phihat}), corresponds to the 
recursion for the continued fraction expansion, and since this map is birational, it is also possible to reverse 
the direction of iteration and extend to all negative indices $n$ (again, this is 
always possible subject 
to the condition that $d_n$ does not vanish for some $n$). This 
immediately leads to a 
J-fraction expression for $Y_0$, that is 
\beq\label{gneg} 
G^\dagger = Y_0 = 
\cfrac{ 2u_0^{-1}d_0}{X+v_{-1} -\cfrac {d_{-1}}{ X+v_{-2}-\cfrac{d_{-2}}{X+v_{-3} -\cdots}} }, 
\eeq 
which corresponds to a power series expansion around $\infty_2$, 
\beq\label{gdags} 
G^\dagger=\sum_{j\geq 1}s_{j-1}^\dagger \, X^{-j} = s_0^\dagger X^{-1}+s_1^\dagger X^{-2}+s_2^\dagger X^{-3}+\ldots .  
\eeq 
This means that the quantities $d_n$, $v_n$ can also be written in terms of ratios of determinants 
when $n$ is negative, but involving the Hankel determinant
$$ 
\Delta_n^\dagger = \det (s_{i+j-2}^\dagger)_{i,j=1,\ldots, n},
$$  
as well as  the associated shifted Hankel determinant $\Delta_n^{\dagger *}$, which is just the 
analogue of (\ref{cas}) built from the coefficients in the series  (\ref{gdags}). 

\begin{thm}\label{hankneg}
For negative indices $n$,  
the quantities $d_n$, $v_n$ that appear under iteration of the J-fraction expansion 
(\ref{gneg}) of the 
function (\ref{pick}), which provide the two components (\ref{topPQ}) 
of the iterates of the birational map (\ref{phi}), are given in 
terms of  Hankel and shifted Hankel determinants  by 
\beq\label{hanknegform} 
d_{1-n} = \frac{\Delta_n^\dagger\Delta_{n-2}^\dagger}{(\Delta_{n-1}^\dagger)^2} \quad\mathrm{for}\,\,n\geq 2, 
\qquad v_{-n} = \frac{\Delta_{n-1}^{\dagger *}}{ \Delta_{n-1}^\dagger}
- \frac{\Delta_{n}^{\dagger *}}{ \Delta_{n}^{\dagger}} \quad \mathrm{for}\,\,n\geq 1, 
\eeq 
where the entries 
are determined recursively by $s_0^\dagger=-u_{-1}/2$ and,  for $j\geq 1$, 
\beq\label{sjdag} 
s_j^\dagger = -\sum_{i=2}^{\rg+1}(k^{(\rg+1-i)} + \pi_0^{(\rg+1-i)}) s_{j-i}^\dagger 
+\frac{d_0}{s_0^\dagger}\left(\sum_{i=0}^{j-2} s_i^\dagger s_{j-2-i}^\dagger 
+\sum_{\ell=3}^{\rg+2}\sum_{i=0}^{j-\ell}\rho_0^{(\rg+2-\ell)} s_i^\dagger s_{j-\ell-i}^\dagger\right) 
+\bar{s}_j^\dagger
\eeq 
where 
$$ 
\bar{s}_j^\dagger= 
\begin{cases}   
s_0^\dagger \rho_{-1}^{(\rg-j)} & \mathrm{for }\,\, 1\leq j\leq \rg, \\ 
0 &\mathrm{for } \,\, j\geq \rg+1. 
\end{cases} 
$$ 
\end{thm} 

\begin{prf} Since $G^\dagger=Y_0=(Y+P_0)/Q_0$, the generating function for the moments $s_j^\dagger$ satisfies the 
quadratic equation 
$P_0G^\dagger -\frac{1}{2}Q_0(G^\dagger)^2 = -\frac{1}{2}Q_{-1}$, analogous to (\ref{genfn}), 
and the recurrence (\ref{sjdag}) 
follows immediately after substituting in the series (\ref{gdags}). Subject to suitable relabelling of indices,  the derivation  
of the formulae 
is the same as in the proof of Theorem \ref{hankt}.   
\end{prf} 

\begin{exa}\label{g2neg} 
{\bf Moments for negative $n$  in genus 2:}
Using  the notation of Example \ref{genus2} once again, 
the recursion (\ref{sjdag}) for $\rg=2$ has initial values  
$$s_0^\dagger=-u_{-1}/2, \qquad  
s_1^\dagger=-s_0^\dagger v_{-1}, \qquad s_2^\dagger = s_0^\dagger(w_{-1}-d_0-f),$$  and 
subsequent coefficients in the series expansion (\ref{gseries}) are 
determined for   $j\geq 3$ by 
\beq\label{g2dagrec} 
s_j^\dagger = -(2d_0+f)s_{j-2}^\dagger -(2d_0e_0+g)s_{j-3}^\dagger 
+(s_0^\dagger)^{-1}d_0\left( 
\sum_{i=0}^{j-2} s_i^\dagger s_{j-2-i}^\dagger -v_0 \sum_{i=0}^{j-3} s_i^\dagger s_{j-3-i}^\dagger
+w_0 \sum_{i=0}^{j-4} s_i^\dagger s_{j-4-i}^\dagger\right) 
. 
\eeq 
In particular, taking the specific curve 
$Y^2=(X^3-5X-1)^2-4(X^2+2X+3)$ 
that was used for illustration in Example \ref{g2hankel},  
with the same function $Y_0$, as before we have 
$f=-5$, $g=u=-1$, 
$u_0=-4$, $d_0=5/4$, $e_0=3/5$, $v_0=-1/2$, $w_0=-3/2$, 
and also $v_{-1}=-1/10$, $w_{-1}=-3/2$,   
which gives 
$$ 
G^\dagger=Y_0 = -\frac{1}{4}\left(\frac{Y+X^3-\frac{5}{2}X+\frac{1}{2}}{X^2+\frac{1}{2}X-\frac{3}{2}}\right)  
=-\frac{5}{8}X^{-1}-\frac{1}{16}X^{-2}-\frac{45}{32}X^{-3}-\frac{25}{64}X^{-4}-\frac{757}{128}X^{-5}-\frac{801}{256}X^{-6}-\cdots, 
$$ 
with $X\to\infty$ and $Y\sim -X^3$,  and the recursion for the coefficients (moments)  in the above 
expansion   is 
$$ 
s_j^\dagger=\frac{5}{2}s_{j-2}^\dagger-\frac{1}{2} s_{j-3}^\dagger 
- 2 \sum_{i=0}^{j-2} s_i^\dagger s_{j-2-i}^\dagger - \sum_{i=0}^{j-3} s_i^\dagger s_{j-3-i}^\dagger
+3 \sum_{i=0}^{j-4} s_i^\dagger s_{j-4-i}^\dagger, \quad j\geq 3
$$  
with initial values $s_0^\dagger=-5/8$, $s_1^\dagger=-1/16$, $s_2^\dagger=45/32$. In this case,  
the sequence of moments $(s_j^\dagger)$ begins 
$$ 
-\frac{5}{2^3},-\frac{1}{2^4},-\frac{45}{2^5},-\frac{25}{2^6},-\frac{757}{2^7},
-\frac{801}{2^8},-\frac{14749}{2^9},-\frac{24361}{2^{10}},-\frac{316037}{2^{11}},-\frac{714865}{2^{12}},\ldots , 
$$ 
yielding the corresponding sequence of Hankel determinants 
\beq\label{g2neghank} 
(\Delta_j^\dagger): \quad 
1,-\frac{5}{2^3},\frac{7}{2^3},-\frac{303}{2^7},\frac{4091}{2^9},-\frac{63805}{2^{10}},\frac{3496637}{2^{12}}, \ldots ,  
\eeq 
which has  alternating signs. 
\end{exa}

\begin{remark} Given the two sets of formulae (\ref{hankform}) and (\ref{hanknegform}), 
it is natural to want to write $d_n$ and $v_n$ in the form  
\beq\label{tauform} 
d_n = \frac{\tau_n\tau_{n-2}}{\tau_{n-1}^2}, \qquad 
v_n = \frac{\tau_{n-1}^*}{\tau_{n-1}} -  \frac{\tau_{n}^*}{\tau_{n}} 
\eeq  
for all $n\in\Z$, for some set of quantities $\tau_n, \tau_n^*$.  However, in general 
one cannot just take $\tau_n=\Delta_n$ for non-negative $n$  and $\tau_n = \Delta_{-n-1}^\dagger$ 
for negative $n$ (and similarly for $\tau_n^*$) , because there will be a mismatch at the values of 
$d_0,d_1$ and $v_0$ which are left unspecified by Theorems \ref{hankt} and \ref{hankneg}. Nevertheless, 
one can make use of the fact that the expressions for $d_n$ and $v_n$ in  (\ref{tauform}) 
are left invariant by the three-parameter group of gauge transformations 
\beq\label{gauge} 
\tau_n \rightarrow ab^n \tau_n , \qquad \tau_n^*\rightarrow ab^n( \tau_n^* + c\tau_n).  
\eeq 
In particular, the choice 
$$ 
\tau_n = 
\begin{cases}   
\Delta_n & \mathrm{for }\,\, n\geq 0, \\ 
(-1)^n \left(\frac{2}{u_0}\right)^{2n+1} \Delta_{-n-1}^\dagger  &\mathrm{for } \,\, n\leq -1
\end{cases} 
$$  
ensures that the values of $d_0$ and $d_1$ match up, and a similar choice can be made 
for $\tau_n^*$ to fix $v_0$. For instance, applying this choice to glue together  
the sequences (\ref{g2hank}) and (\ref{g2neghank}) in a consistent fashion 
yields the doubly infinite  sequence 
\beq\label{doubtau} 
(\tau_n): \quad \ldots, 
562196701, 6993274, 127610, 4091, 303, 28, 5,2,1,1,2,6,31,319, 5810, 147719, 
\ldots .  
\eeq  
\end{remark} 

\section{The Somos connection}\label{somoscon} 
\setcounter{equation}{0}

In this section, we explain how Somos sequences naturally arise from the continued fraction expansion, 
as quadratic relations for the Hankel determinants $\Delta_n$. This is most straightforward 
to describe in the genus one  case, as it follows from the fact that, 
for a fixed value of the first integral $H=-uv$ given by (\ref{g1ham}), 
each orbit of (\ref{g1map}) coincides with an orbit of a symmetric QRT map, and, as was already noted in 
\cite{qrt} in an example related to the hard hexagon model, the bilinear form of the latter is 
precisely (\ref{s4recu}). (For a detailed discussion of normal forms of QRT maps restricted 
to fixed invariant curves, see \cite{ir,ir2}.) 

\begin{propn}\label{s4qrt}
For a fixed value of the first integral $H=-uv$, on each orbit of the map (\ref{g1map}) the 
quantity $d_n$ satisfies the second order recurrence 
\beq\label{qrtmap}
d_{n+1}d_{n-1}=\frac{\al\, d_n +\be}{d_n^2},  
\eeq  
with coefficients $\al = u^2$,  $\be = u^2(v^2+f)$.
\end{propn} 
\begin{prf}
Putting $X=v_n$ into (\ref{g1}) and then applying the first equation in (\ref{g1rels}) yields 
$$ 
u(v_n-v) = d_n (d_n+v_n^2+f) =-d_nd_{n+1}. 
$$
Then putting $X=v$  into (\ref{g1}) and using the above result gives 
$$ 
\begin{array}{rcl} 
0 & = & d_n \Big( v^2+f - (v-v_{n-1})(v-v_{n})+d_n \Big) \\ 
& = &  d_n \Big( v^2+f - u^{-2}d_{n+1}d_n^2 d_{n-1}+d_n \Big),  
 \end{array} 
$$ 
so that (\ref{qrtmap}), which is an example of 
a symmetric QRT map \cite{qrt},  follows immediately. 
\end{prf} 

The connection with Somos-4 is almost immediate, since if $d_n =\tau_n\tau_{n-2}/(\tau_{n-1})^2$, 
then $\tau_n$ satisfies (\ref{s4recu}) whenever   
$d_n$ satisfies (\ref{qrtmap}). So in particular, by Theorem \ref{hankt}, the Hankel
determinants for $\rg=1$ satisfy a Somos-4 relation, and since the latter is invariant under 
the first of the  
gauge transformations (\ref{gauge}), any Somos-4 sequence can be expressed in terms of Hankel   
determinants. More precisely, starting from any Somos-4 sequence, one can always  make a gauge 
transformation to a sequence with $\tau_0=1$, and then use the coefficients 
$\al,\be$  and the 
other initial conditions $\tau_1,\tau_2,\tau_3$ to specify the values of $s_0=\Delta_1=\tau_1$ 
and $v_0,v_1,d_1,f$, so that the values of $s_1$ and the 
coefficients in (\ref{elliprec}) are fixed. (In fact, since the gauge transformation 
involves two parameters $a,b$, there is also the freedom to fix 
$\tau_1=1$, which corresponds to taking $s_0=1$.)  Thus we arrive at the following result. 

\begin{thm}\label{s4h} 
In the case $\rg=1$, the Hankel determinants (\ref{hankel}) 
with moments defined recursively by (\ref{elliprec}) 
satisfy the Somos-4 recurrence 
\beq\label{s4hank} 
\Delta_{n+4}\Delta_n = \al \, \Delta_{n+3}\Delta_{n+1} + \be\, \Delta_{n+2}^2, 
\eeq  
with 
\beq \label{abvals}  
\al = u^2, \qquad \be = u^2(v^2+f). 
\eeq 
Moreover, 
every 
solution of the Somos-4  recurrence (\ref{s4recu}) can be written in the form 
$$ 
\tau_n = 
\begin{cases}   
\hat{a}\hat{b}^n \Delta_n & \mathrm{for }\,\, n\geq 0, \\ 
a^\dagger (b^\dagger)^n \Delta_{-n-1}^\dagger  &\mathrm{for } \,\, n\leq -1, 
\end{cases} 
$$  
where $\Delta_{n}^\dagger$ is constructed from moments that satisfy 
(\ref{sjdag}) with $\rg=1$, 
for suitable constants $\hat{a},\hat{b},a^\dagger,b^\dagger$. 
\end{thm}  

There is an apparent mismatch between the Hankel determinants in (\ref{dhank}), which 
were shown by Xin to yield the terms of the Somos-4 sequence (\ref{s4seq}), and those 
in (\ref{haseq}) above. We now explain the relation between these two sets of 
Hankel determinant formulae, and see how the results of \cite{chx}  
are a consequence of the continued fraction expansion for $\rg=1$. 

\begin{thm}\label{dhankthm} 
For $n\geq 2$, the  quantity $d_n$ that satisfies  (\ref{g1map}) is given by 
\beq\label{ddhank}  
d_n = \frac{D_nD_{n-2}}{D_{n-1}^2} , 
\eeq 
in terms of the Hankel determinant (\ref{origh}) defined in terms of moments $\tilde{s}_j$ 
that satisfy the recursion (\ref{g1altrec}) for $j\geq 2$, with 
$$\tilde{s}_0=u_1/2, \quad \tilde{s}_1=-\tilde{s}_0(v_0+v_1), 
\quad \tilde{\al}=-2v_0, \quad 
\tilde{\be}=d_0-d_1, \quad 
\tilde{\gam}=\frac{d_1}{\tilde{s}_0}.
$$  
Moreover, the sequence $(D_n)_{n\geq 0}$ is identical to 
the sequence of Hankel determinants $(\Delta_n)_{n\geq 0}$ 
with moments satisfying (\ref{elliprec}), hence 
satisfies the Somos-4 recurrence 
(\ref{s4hank}) with coefficients   
as in (\ref{abvals}). 
\end{thm} 
\begin{prf}By replacing $X$ with the shifted variable $\tilde{X}=X-v_0$, and 
letting $\tilde{G}(\tilde{X}) = G(X)$, we obtain the J-fraction  
\beq\label{gtcla} 
\tilde{G}
=\cfrac{ -2\tilde{u}_0^{-1}\tilde{d}_1}{\tilde{X}+\tilde{v}_1 -\cfrac {\tilde{d}_2}{ \tilde{X}+\tilde{v}_2-\cfrac{\tilde{d}_3}{\tilde{X}+\tilde{v}_3 -\cdots}} }= \sum_{j\geq 1} \tilde{s}_{j-1}\tilde{X}^{-j},   
\eeq
where, from (\ref{genfn}) in the case $\rg=1$, the generating function $\tilde{G}$ satisfies 
\beq\label{gtquad} 
\Big(\tilde{X}^2+2v_0\tilde{X}+(v_0^2+2d_1+f)\Big) \tilde{G} +\frac{u_0}{2} \tilde{X} \tilde{G}^2 
=\frac{u_1}{2}(\tilde{X}+v_0-v_1). 
\eeq 
In the first line of the continued fraction (\ref{gtcla}), we are at liberty to choose 
$\tilde{u}_0=u_0$, which implies  that $\tilde{d}_1=d_1$, and then in each 
subsequent line we have $\tilde{v}_n = v_n +v_0$ for $n\geq 1$, and also $\tilde{d}_n=d_n$ 
for $n\geq 2$. Hence, by the same argument as in the proof of Theorem \ref{hankt},  
the formula (\ref{ddhank}) holds for $n\geq 2$. The moments $\tilde{s}_j$ are obtained 
from the series expansion of $\tilde{G}$ in powers of $\tilde{X}$, with the 
leading order term (order $\tilde{X}$) in (\ref{gtquad}) giving $\tilde{s}_0=s_0 = u_1/2$, and the next to leading 
order term (order $\tilde{X}^0$) giving 
$\tilde{s}_1 +2v_0\tilde{s}_0  =\frac{u_1}{2}(v_0-v_1)$,  
while at order $\tilde{X}^{-j+1}$ for $j\geq 2$, 
upon noting that $\tilde{\be} = -v_0^2 -2d_1-f = d_0-d_1$ from (\ref{g1map}), we find the recursion relation (\ref{g1altrec}) with the stated values of $\tilde{\al},\tilde{\be},\tilde{\gam}$.  
By convention we have $D_0 = \Delta_0=1$, and also $D_1=\tilde{s}_0=s_0=\Delta_1$, 
and then it follows by induction from (\ref{ddhank}) that $D_n=\Delta_n$ for all $n\geq 0$. 
\end{prf} 

\begin{remark} In order to see how Xin's result \cite{xin} follows from the above, it is sufficient 
to note that the quartic curve (\ref{s4quart}) is isomorphic to (\ref{37a}), via 
the birational equivalence $Y=x^{-2}(1-2\eta)$, $X=\tilde{X}-1=x^{-1}-1$, so that the 
expansion (\ref{gtilde}) in powers of $x$ is equivalent to an expansion in powers of $\tilde{X}^{-1}$. Also, 
by setting $\eta=(1-y)/2$,  the Weierstrass cubic $y^2=4x^3-4x+1$  derived from analytic formulae in \cite{honeblms} is seen to be isomorphic to the curve (\ref{37a}); over $\Q$, this is known as 
37a1, the elliptic curve  of 
minimal conductor  with positive rank. (See {\tt www.lmfdb.org/EllipticCurve/Q/37/a/1} in the online 
database of L-functions, modular forms, and related objects.)  
\end{remark} 

The analogue of Theorem \ref{s4h} in genus two is more difficult to state explicitly, due to 
the size of the expressions for the coefficients, and at present we are  only able to prove it 
with the use of computer algebra.  

\begin{thm}\label{g2h} 
In the case $\rg=2$, the Hankel determinants (\ref{hankel}) 
with moments defined recursively by (\ref{g2rec}) 
satisfy a Somos-8 recurrence of general type, that is   
\beq\label{g2s8} 
\al_1\, \Delta_{n+8}\Delta_n + \al_2 \, \Delta_{n+7}\Delta_{n+1} + \al_3\, \Delta_{n+6}\Delta_{n+2} 
+ \al_4\, \Delta_{n+5}\Delta_{n+3}+ \al_5\, \Delta_{n+4}^2=0, 
\eeq  
where the coefficients $\al_1,\ldots,\al_5$ are certain first integrals 
of the map (\ref{g2pair}). 
\end{thm}
\begin{prf} 
The recurrence (\ref{g2s8}) is equivalent to a relation for the iterates  $d_n$, that is 
$$ 
\al_1 d_{n+8}d_{n+7}^2d_{n+6}^3d_{n+5}^4d_{n+4}^3d_{n+3}^2d_{n+2}
+\al_2 d_{n+7}d_{n+6}^2d_{n+5}^3d_{n+4}^2d_{n+3}
+\al_3 d_{n+6}d_{n+5}^2d_{n+4} 
+\al_4 d_{n+5}+\al_5 =0,  
$$  
along an orbit of the 4D map (\ref{g2pair}). 
Equivalently, writing a solution of this map in the form 
(\ref{tauform}), it means that $\tau_n$ should satisfy the Somos-8  
recurrence (\ref{g2s8}) with some  coefficients $\al_j$ that are constant 
on each orbit. This requires the vanishing of a $5\times 5$ determinant, namely   
\beq \label{bigdet} 
\left| \begin{array}{ccccc} 
\tau_{n+4}\tau_{n-4} & \tau_{n+3} \tau_{n-3} &  \tau_{n+2} \tau_{n-2} &\tau_{n+1} \tau_{n-1} &
\tau_{n}^2 \\ 
\tau_{n+5}\tau_{n-3} & \cdots & \cdots & \cdots & \tau_{n+1}^2 \\ 
\vdots  &&&&\vdots \\ 
\tau_{n+8}\tau_{n} & \cdots & \cdots & \cdots & \tau_{n+4}^2
\end{array}\right| =0
\eeq  
for all $n$, and also that the ratios of certain $4\times 4$ minors should be independent of $n$. 
In particular, denoting by $\hat{D}_{j,n}$ the minor formed from the first four rows in (\ref{bigdet}) with the 
$j$th column removed, so that  
$$ 
\hat{D}_{1,n}=
\left| \begin{array}{cccc} 
 \tau_{n+3} \tau_{n-3} &  \tau_{n+2} \tau_{n-2} &\tau_{n+1} \tau_{n-1} &
\tau_{n}^2 \\ 
\tau_{n+4}\tau_{n-2} &  \cdots & \cdots & \tau_{n+1}^2 \\ 
\tau_{n+5}\tau_{n-1} & \cdots &\cdots & \tau_{n+2}^2\\ 
\tau_{n+6}\tau_{n} & \cdots & \cdots &  \tau_{n+3}^2
\end{array}\right|,   
$$ 
and so on, the existence of the relation (\ref{g2s8}) for $\tau_n$  is equivalent 
to the  requirement that the ratios 
$$ 
\frac{\al_j}{\al_1} = \hat{D}_{j,n} /\hat{D}_{1,n}, \qquad j=2,3,4,5 
$$ 
are independent of $n$ (noting 
the possibility of a vanishing denominator, in the case that $\al_1=0$). 
So in order to verify the statement, it is sufficient to check that 
\beq\label{detd} 
\hat{D}_{j,n}\hat{D}_{1,n+1}-\hat{D}_{j,n+1}\hat{D}_{1,n}=0, 
\eeq  
holds for each $j=2,3,4,5$, and for all $n$. In fact, for each $j$, it is enough to check that this holds 
for a single shift $n\to n+1$ with arbitrary initial data in  the map (\ref{g2pair}), which 
guarantees that each of the ratios $\al_j/\al_1$ is a first integral, and then 
it automatically holds for all $n$. Even with the help of  computer algebra, this is 
not a completely straightforward task, and  in order to do it as efficiently as possible it is convenient 
to  make a gauge transformation (\ref{gauge}) to fix $s_0=1$, and then note that 
there is a one-to-one correspondence between two sets of seven parameters:
the four initial values $d_0,d_1,v_0,v_1$ and three parameters $f,g,u$ 
needed to iterate the map        (\ref{g2pair}), and the two  initial values $s_1,s_2$ 
and five coefficients $\tilde{\al},\tilde{\be},\tilde{\gam},\tilde{\delta},\tilde{\epsilon}$ that 
specify the genus two  
recursion  
(\ref{g2rec}) in the form  
$$ 
s_j = \tilde{\al} s_{j-2}+\tilde{\be} s_{j-3} 
+ 
\tilde{\gam} \sum_{i=0}^{j-2} s_i s_{j-2-i} + \tilde{\delta}\sum_{i=0}^{j-3} s_i s_{j-3-i}
+ \tilde{\epsilon}\sum_{i=0}^{j-4} s_i s_{j-4-i} 
. 
$$
Now, using the above recursion, one can calculate the sequence $(s_j)_{j\geq 0}$, and 
then compute the  Hankel determinants $\tau_n=\Delta_n$ for $n=0,1,2,\ldots$, 
which are polynomials in $\Z[s_1,s_2,\tilde{\al},\tilde{\be},\tilde{\gam},\tilde{\delta},\tilde{\epsilon}]$, 
but this rapidly becomes very 
computationally intensive as $n$ increases. More efficient is to rewrite   the map   (\ref{g2pair}) as an equivalent  
pair of coupled recurrence relations of  degree 6  
for $\tau_n, \tau_n^*$, which are of overall  order 7. To iterate the latter, one needs 7 initial values (four adjacent 
$\tau_n$ and three adjacent $\tau_n^*$), and it is convenient to take 
$\tau_0=\Delta_1=1$, $\tau_1=\Delta_0=s_0=1$ but also $\tau_{-1}=d_1=\tilde{\gam}$ 
and 
$\tau_{-2} = \tilde{\gam}^2 d_0 = \tilde{\al}\tilde{\gam}+ 
\tilde{\gam}^3 - \tilde{\delta}^2+\tilde{\gam}\tilde{\epsilon}
$, 
together with $\tau_0^*=\Delta_0^*=0$, $\tau_1^*=\Delta_1^*=s_1$, as 
well as $\tau_{-1}^* = d_1v_0 =-\tilde{\delta}$, 
noting that $f=-\tilde{\al} -2\tilde{\gam}$, $g=-\tilde{\be}-2\tilde{\gam}s_1-2\tilde{\delta}$, 
$u=-\tilde{\gam}s_2-\tilde{\delta}s_1-\tilde{\epsilon}$. The verification of 
(\ref{detd}) requires 13 adjacent values of $\tau_n$, but due to the size of the 
expressions involved it is best to compute only up to $\tau_6$ using five forward steps of the coupled recurrence 
for   $ \tau_n, \tau_n^*$, and then apply  this recurrence in reverse, making four backward steps to go back as 
far as $\tau_{-6}$, so that the adjacent values $\tau_{-6},\tau_{-5},\ldots, \tau_5,\tau_6$ 
are obtained as explicit polynomials in  $\Z[s_1,s_2,\tilde{\al},\tilde{\be},\tilde{\gam},\tilde{\delta},\tilde{\epsilon}]$. 
This means that the minors $\hat{D}_{j,-2}$ and $\hat{D}_{j,-1}$ can be computed explicitly, 
which allows (\ref{detd}) to be checked directly when $n=-2$. The formulae for the first 
integrals $\al_j/\al_1$, $j=2,3,4,5$ as rational functions of $s_1,s_2,\tilde{\al},\tilde{\be},\tilde{\gam},\tilde{\delta},\tilde{\epsilon}$ 
are so large that they are difficult to display even on a computer screen, but 
if $\al_1$ is regarded as the first integral that is the lowest common denominator 
of these four quantities, then we arrive at the Somos-8 relation (\ref{g2s8}). 
\end{prf} 

\begin{exa} 
The doubly infinite sequence (\ref{doubtau}), which extends (\ref{g2hank}), satisfies the Somos-8 relation
$$
7\tau_{n+8}\tau_{n} +137\tau_{n+7}\tau_{n+1} +2504\tau_{n+6}\tau_{n+2} 
-43424\tau_{n+5}\tau_{n+3}
-26959\tau_{n+4}^2=0.
$$
\end{exa} 

\begin{remark}
As already noted, there is the possibility of a vanishing denominator $\al_1=0$ in 
the ratios $\al_j/\al_1$, $j=2,3,4,5$. Given 
that the map (\ref{g2pair}) only has two 
independent first integrals, which can be specified by $H_1=-uv,H_2=uw$ as  in  
(\ref{nrel3}) and (\ref{nrel4}),   it follows that these four ratios are rational functions 
of $f,g,u,H_1,H_2$, with coefficients in $\Q$, so that $\al_1$ can be fixed as the polynomial 
in $\Z[f,g,u,H_1,H_2]$ that is the lowest common denominator 
of these  four rational functions.  Thus it can happen that $\al_1 =0$ for certain combinations 
of $f,g,u,H_1,H_2$, in which case $\Delta_n$ (or $\tau_n$) satisfies a Somos-6 relation, rather than a Somos-8. 
Numerical experiments suggest that $u|\al_1(f,g,u,H_1,H_2)$,  consistent with a result of 
van der Poorten, who showed that there is a Somos-6 relation in the special case $u=0$ \cite{vdp2,vdp3}. 
\end{remark} 

Before  concluding this section, we state a conjecture which is the genus two  analogue of Theorem \ref{s4sigma}. 

\begin{conje}\label{g2sigma} 
When $\rg=2$, the Hankel determinants (\ref{hankel}) 
with moments defined  by (\ref{g2rec})  
are given by 
\beq\label{g2sig} 
\Delta_n=\hat{a}\hat{b}^n\, \frac{\si ({\bf z}_0+n{\bf z})}{\si ({\bf z})^{n^2}}, 
\eeq  
in terms of  the 
genus two Kleinian sigma function $\si ({\bf z})=\si ({\bf z};\tilde{c}_0,\tilde{c}_1,\tilde{c}_2,\tilde{c}_3)$ 
associated with a quintic curve $\tilde{{\cal C}}$ that is isomorphic to 
the sextic 
${\cal C}$ in (\ref{csextic}), given by  
$\tilde{{\cal C}}: \, y^2=4x^5+\sum_{j=0}^3\tilde{c}_jx^j$  with period lattice $\Lambda$, for 
${\bf z}_0$, ${\bf z}\in\C^2\bmod\Lambda$ 
with ${\bf z} =2 \int_\infty^{\tilde{P}_2}(\frac{\rd x}{y}, \frac{x\rd x}{y} )^T$, where 
$\infty$ is the unique point at infinity on  $\tilde{{\cal C}}$, $\tilde{P}_2\in \tilde{{\cal C}}$ is the point 
corresponding to  $\infty_2\in {\cal C}$, and $\hat{a},\hat{b}$ are certain non-zero constants.  
\end{conje} 

A proof of the above result would follow from an analytic solution for the iterates of the map 
 (\ref{g2pair}), which we propose to consider elsewhere. However, to see why this result is 
plausible we note that since the class of the divisor  $\tilde{P}_2-\tilde{P}_1\sim_{\mathrm{lin}} 2 (\tilde{P}_2-\infty)$ 
 in the Jacobian of $\tilde{{\cal C}}$ corresponds to that of the divisor 
$\infty_2 -\infty_1$ on  ${\cal C}$, each shift $n\to n+1$ increases the argument of 
the numerator in (\ref{g2sig}) by ${\bf z}$, which is consistent with (\ref{lineq}). 
Furthermore, if the formula (\ref{g2sig}) is correct then, by essentially the same analytical calculations as  
those in \cite{beh, hones6}, it follows immediately that $\Delta_n$ satisfies the Somos-8 
relation (\ref{g2s8}), or a Somos-6 relation when a certain constraint on ${\bf z}$ holds.

In the  higher genus case, we further conjecture that there should be an analytical formula analogous to 
(\ref{g2sig}). Equivalently, there should be an expression in terms of the Riemann theta function 
associated with the Jacobian of the curve (\ref{curve}), of the 
form $\Delta_n = \hat{a}\hat{b}^n\hat{c}^{n^2} \Theta ({\bf v}_0+n{\bf v})$ 
for some $\hat{a},\hat{b},\hat{c}\in\C^*$ and ${\bf v}_0,{\bf v}\in \C^{\rg}$. If this 
expression holds, then by counting 
the dimension of the vector space of quasiperiodic functions of weight 2 with respect to the period lattice 
(see \cite{tata1}), it follows that $\Delta_n$ satisfies a Somos-$k$ relation for some $k\leq 2^{\rg+1}$. 
We have verified this in numerical examples for $\rg=3,4$.

\section{Poisson structure and integrability} \label{poisson} 
\setcounter{equation}{0}
In this section we slightly change our notation, and consider a modified family of  Lax matrices, given by 
\beq\label{newlax}  
{\cal L}(\ze)= \left(\begin{array}{cc} \rP(\ze) & \rR(\ze) \\ 
\rQ(\ze) & -\rP(\ze) \end{array}\right),  
\eeq 
where $\rP$ and $\rR$ are monic polynomials of degrees $\rg+1$ and $\rg$ in $\ze$, 
respectively, and $\rQ$ is a polynomial of degree $\rg$ with non-constant leading coefficient, 
which we write as 
$$ 
\rQ = -4d_0\ze^\rg+ O(\ze^{\rg-1}). 
$$  
The set of all such matrices forms an affine space of dimension $3\rg+2$, with 
coordinates given by the non-trivial coefficients of $\rP$, $\rQ$  and $\rR$. 
We will endow this space with a particular Poisson structure of rank $2\rg$, and show how this 
leads to the construction of an associated set of Hamiltonian vector fields that 
are completely integrable in the Liouville sense. Then we will present a compatible 
discrete integrable system on the same phase space, and show that this is equivalent to 
the iteration of the recursion for the continued fraction expansion of the 
hyperelliptic function considered 
previously.      

The Poisson brackets between the entries of $\cal L$ are specified by  
\beq\label{pbs1} 
\{\, \rP(\ze), \rP(\eta)\,\}  =  0 \, = \, \{\, \rR(\ze), \rR(\eta)\,\}, 
\eeq 
 \beq\label{pbs2} 
\{\, \rP(\ze), \rQ(\eta)\,\} =  2\left(\frac{\rQ(\ze)-\rQ(\eta) }{\ze-\eta}\right)  , 
\quad \{\, \rP(\ze), \rR(\eta)\,\}  = - 2\left(\frac{\rR(\ze)-\rR(\eta) }{\ze-\eta}\right)  ,
\eeq 
\beq\label{pbs3} 
\{\, \rQ(\ze), \rQ(\eta)\,\} =  -4\Big(\rQ(\ze)-\rQ(\eta)\Big)  , 
\quad
\{\, \rQ(\ze), \rR(\eta)\,\} =  4\left(\frac{\rP(\ze)-\rP(\eta) }{\ze-\eta}\right) - 4 \rR(\eta).   
\eeq In terms of the coefficients 
appearing in $\cal L$, this is a linear bracket, since the right-hand sides are linear in 
$\rP$, $\rQ$, $\rR$.  
In order for this to define a Poisson bracket, it must  satisfy the Jacobi identity, and 
although this can be verified directly 
this requires many  
tedious calculations; we set this question aside for now, and a simpler argument 
will be presented in due course. 

To begin with, we 
consider the function 
\beq\label{Fdef}
F(\ze) = -\det {\cal L} (\ze)=\rP(\ze)^2 +\rQ(\ze)\rR(\ze). 
\eeq 
Using the above bracket relations, we find that 
\beq\label{pf} 
\{\, \rP(\ze), F(\eta)\,\}  =2\left(\frac{\rQ(\ze)\rR(\eta)-\rQ(\eta)\rR(\ze)}{\ze-\eta}\right), 
\eeq 
\beq\label{qf} 
\{\, \rQ(\ze), F(\eta)\,\}  =4\left(\frac{\rP(\ze)\rQ(\eta)-\rP(\eta)\rQ(\ze)}{\ze-\eta}\right) 
-4\rQ(\ze)\rR(\eta), 
\eeq 
\beq\label{rf}
\{\, \rR(\ze), F(\eta)\,\}  =4\left(\frac{\rP(\eta)\rR(\zeta)-\rP(\zeta)\rR(\eta)}{\ze-\eta}\right) 
+4\rR(\ze)\rR(\eta). 
\eeq 
It is straightforward to check that the right-hand sides of the above expressions 
are polynomials in $\eta$ of degrees $\rg-2$, $\rg-1$, $\rg-1$ respectively. Thus, 
if we expand 
\beq\label{fexpc}
F(\eta)=\eta^{2\rg+2}+\sum_{j=0}^{2\rg+1}c_j\ze^j, 
\eeq 
then these expressions imply that  
$$ 
\{ \, \cdot \, , c_j\,\}=0 \quad \mathrm{for} \quad j=\rg,\rg+1,\ldots,2\rg+1, 
$$ or in other words the leading $\rg+2$ non-trivial coefficients of $F$ 
are Casimirs. 

In fact $c_\rg, \ldots, c_{2\rg+1}$ provide the  full set of  Casimirs, and the symplectic leaves 
have dimension $2\rg$. To see this, 
factorize $\rR$ as 
\beq\label{xdef} \rR(\ze)=\prod_{i=1}^\rg (\ze-x_i), \eeq  
and set 
\beq\label{ydef} 
y_i=\rP(x_i), \qquad i=1,\ldots,\rg. 
\eeq 
Then 
$$ 
\rP(\ze)=\ze^{\rg+1}+\frac{1}{2}c_{2\rg+1}\ze^\rg +\sum_{j=0}^{\rg-1}\pi_j\ze^j, 
$$ 
where, from (\ref{ydef}), the $\rg$ coefficients $\pi_j$ can be found explicitly as functions of 
the $x_i$, the $y_i$ and the Casimir    $c_{2\rg+1}$ by solving a linear system. 
The two relations in (\ref{pbs1}) then imply that 
$$ 
\{ \, x_i, x_j \, \} =0 = \{\,y_i, y_j\,\} 
$$ 
for all $i,j$, 
while from  the second bracket in (\ref{pbs2}) it follows that 
$$ 
\{ \, y_i, x_j\, \} =2\delta_{ij}, 
$$ 
so that (up to scaling) the pairs $(x_i,y_i)_{i=1,\ldots,\rg}$ 
provide a set of canonical coordinates on a symplectic manifold of 
dimension $2\rg$, and by (\ref{Fdef}), for fixed 
values of the coefficients $c_j$ in (\ref{fexpc}), each pair $(\ze,\mu)=(x_i,y_i)$ is a point  on the spectral curve 
\beq\label{muspec} 
\mu^2 = F(\zeta).
\eeq  
Then  by using (\ref{xdef}), (\ref{ydef}) and the leading terms   
in (\ref{Fdef}) up to and including order $\ze^g$, all of the coefficients in $\cal L$ can be expressed as 
functions of the $x_i$, $y_i$ and the Casimirs, so the whole Poisson algebra is expressed in terms of 
these coordinates. 

The latter argument begs the question of whether the Poisson brackets for the entries 
of $\cal L$ satisfy the Jacobi identity in the first place, but this is easily seen by reversing the direction 
of the preceding argument. Indeed,  starting with the canonical bracket between the $x_i$ and $y_j$, one  
extends it with the Casimirs $c_\rg, \ldots, c_{2\rg+1}$ to obtain a Poisson algebra of dimension 
$3\rg+2$, where the entries of $\cal L$ are defined in terms of these 
canonical coordinates as above, 
and by construction they satisfy the linear bracket relations given before.  
Hence the Jacobi identity is trivially satisfied.    

We now consider a family of vector fields on the space of Lax matrices, defined 
by the flow 
\beq\label{fflow} 
\frac{\rd}{\rd t}  {\cal L}(\ze) = \{ \, {\cal L}(\ze),F(  \eta)\,\}. 
\eeq 
From the bracket relations (\ref{pf}), (\ref{qf}), (\ref{rf}), it can be verified directly that this can be written in 
the form of a Lax equation, that 
is 
$$ 
\frac{\rd}{\rd t}  {\cal L}(\ze) = [{\cal P}(\ze,\eta), {\cal L}(\ze)],  
$$
with the matrix 
\beq\label{pmat}
{\cal P}(\ze,\eta) = \frac{2}{\ze-\eta}\left(\begin{array}{cc} 
 \rP (\eta) +({\ze-\eta})\rR (\eta) & \rR (\eta) \\ 
 \rQ (\eta) &  -\rP (\eta) -({\ze-\eta}) \rR (\eta) 
 \end{array} \right). 
\eeq
Moreover, all the flows in this family commute with one another, because 
the same bracket relations imply that 
\beq\label{fbf} 
\{\, F(\ze), F(\eta)\,\}  =2\rP(\ze) \{\, \rP(\ze), F(\eta)\,\} + \rR(\ze) \{\, \rQ(\ze), F(\eta)\,\} 
+ \rQ(\ze) \{\, \rR(\ze), F(\eta)\,\}  
=0. 
\eeq 
Thus it follows that all of the coefficients $c_j$ in (\ref{fexpc}) are in involution, and in order to 
get non-trivial flows we can take the non-Casimir functions (Hamiltonians)  
\beq\label{hams} 
H_j = c_j =\mathrm{res}\, \frac{F(\eta)}{\eta^{j+1}} \Big\vert_{\eta=0}, \qquad j=0,1,\ldots, \rg-1. 
\eeq 
Then since $\{\, H_j,H_k\,\}=0$ for all $j,k$, this gives $\rg$ commuting flows which can 
be written in Lax form as 
$$ 
\frac{\rd}{\rd t_j}  {\cal L}(\ze) = \{ \, {\cal L}(\ze),H_j\,\} = [{\cal P}^{(j)}(\ze), {\cal L}(\ze)], 
$$ 
where, for $j=0,\ldots,\rg-1$, the matrix ${\cal P}^{(j)}$ is defined from (\ref{pmat}) by 
$$
{\cal P}^{(j)}(\ze) = \mathrm{res}\, \frac{{\cal P}(\ze,\eta)}{\eta^{j+1}}\Big|_{\eta=0}. 
$$ 
Hence we have integrability in the sense of Liouville \cite{ar}. 

\begin{thm}\label{lint} 
The Hamiltonians (\ref{hams}) define a completely integrable system on the space of Lax matrices (\ref{newlax}). 
\end{thm}

In fact, 
there is more that one can say: for fixed values of $c_j$  in (\ref{fexpc}), the set 
of triples of polynomials $\rP(\ze),\rQ(\ze),\rR(\ze)$ of the specified form that satisfy (\ref{Fdef}) 
for a fixed set of coefficients $c_j$ in (\ref{fexpc}) is an 
affine algebraic variety of dimension $\rg$, that is canonically isomorphic to the affine Jacobian 
of the corresponding hyperelliptic spectral curve (\ref{muspec}), by associating each such 
triple with the degree $\rg$ divisor 
$$ 
D=(x_1,y_1) + \cdots + (x_\rg,y_\rg ) 
$$ 
defined by (\ref{xdef}) and (\ref{ydef}). 
As is explained in \cite{tata2}, 
the analogous construction of the Jacobian variety in the case of odd hyperelliptic curves 
goes back to Jacobi,  and arises in the context of finite gap 
solutions of the Korteweg--deVries equation.  The Poisson brackets and first integrals $H_j$ are all algebraic - in actual fact, they are  
given by polynomial functions of the coefficients of the polynomials $\rP(\ze),\rQ(\ze),\rR(\ze)$ - and over $\C$ the 
generic common level set of these first integrals is an affine part of a complex algebraic torus, with the 
Hamiltonian flows being linear on the torus,  so this is what 
is known as an algebraic completely integrable system (see \cite{vanhaecke} and references for details). 

We now proceed to describe how the Lax matrices (\ref{newlax}) are related to the discrete
Lax pair 
and nonlinear map introduced in section 3. 

First of all, observe that 
completing the square in (\ref{fexpc}) means that 
$$
F=A^2+\tilde{c}_\rg \ze^\rg +\tilde{c}_{\rg-1}\ze^{\rg-1}\cdots +\tilde{c}_0, \quad 
A(\ze)=\ze^{\rg+1} 
+\tilde{c}_{2\rg+1}\ze^\rg+\cdots+\tilde{c}_{\rg+1},
$$
where $\tilde{c}_{2\rg+1}=\frac{1}{2}c_{2\rg+1}$ and all 
the coefficients $\tilde{c}_{2\rg+2-j}$ for $j=1,\ldots,\rg+2$ are polynomial functions 
of the $\rg +2$ Casimirs $\tilde{c}_{2\rg+2-j}$ for the same range of $j$, and there is a bijection between 
these two sets of Casimir functions. In particular, all of the coefficients of $A$ are Casimirs, and 
from (\ref{Fdef}) we may write 
$$ 
\rP(\ze) =A(\ze)+\sum_{j=0}^{\rg-1} \pi_0^{(j)}\ze^j, 
$$ 
and take the $\rg$ non-trivial coefficients of $\rR(\ze)$ and the quantities $\pi_0^{(j)}$ for 
$j=0,\ldots,\rg-1$ as coordinates on each symplectic leaf. We shall see shortly that the latter 
is consistent with the notation in (\ref{pqexp}), but before  getting to this we must restrict 
to a particular set of symplectic leaves, by fixing the value of the top Casimir to be zero, i.e. $ c_{2\rg+1}=0$, 
so that 
$\tilde{c}_{2\rg+1}=0$, and the coefficients of $F$, $A$ and $\rP$ at next to leading order 
are zero; this slightly simplifies the formulae for the nonlinear map, and  agrees with our previous conventions for the continued fraction expansion 
(but if necessary the case of non-zero $ c_{2\rg+1}$ 
can always be obtained by making a shift in the spectral parameter $\ze$). 

Next, in order to obtain (\ref{newlax}), we wish to remove the 
multipliers $u_n$ that appear in the off-diagonal terms of (\ref{lax}), since although they 
provide an arbitrary choice of scale in the continued fraction, they 
do not behave well from the point of view of the Poisson structure. 
Thus we  consider diagonal gauge transformations 
$$
\mathbf{\Psi}_n\mapsto \mathbf{\Phi}_n = {\bf G}_n \mathbf{\Psi}_n, 
\qquad 
{\bf G}_n=\left(\begin{array}{cc} \la_n^{(1)} & 0 \\ 0 & \la_n^{(2)}\end{array}\right) 
$$ 
applied to the eigenvector in (\ref{eigval}).  
These have the effect of changing the prefactors in the off-diagonal entries of    
 ${\bf L}_n$ while leaving the diagonal terms the same.  
Hence, for a suitable choice of $\la_n^{(1)},\la_n^{(2)}$, upon setting $n=0$ we obtain 
$$ 
{\cal L}(\ze ) = {\bf G}_0  {\bf L}_0(\ze ) {\bf G}_0^{-1}
$$ 
with ${\cal L}(\ze )$ being given by (\ref{newlax}), and from (\ref{mshift}) we 
find 
$$ 
{\cal M}(\ze ) =  {\bf G}_0  {\bf M}_0(\ze ) {\bf G}_1^{-1}
= \left(\begin{array}{cc} \ze +v_0 &  1/2 \\ -2d_0 & 0 \end{array}\right),  
$$
where
\beq\label{dvdef} 
d_0 = -\frac{1}{4} \mathrm{res}\, \frac{\rQ(\ze)}{\ze^{\rg+1}}\Big|_{\ze=0}, 
\quad 
v_0= \frac{1}{4d_0} \mathrm{res}\, \frac{\rQ(\ze)}{\ze^{\rg}}\Big|_{\ze=0}. 
\eeq 
Comparing the effect of the gauge transformation with the 
notation used in section 3, it is apparent that  $\rP(\ze)=P_0(\ze)$, 
$\rQ(\ze) = u_{-1}Q_0(\ze)$ and $\rR(\ze)=Q_{-1}(\ze)/u_{-1}$.

Finally, to rewrite the nonlinear map in terms of the   new Lax matrices, we 
set 
$$ 
\tilde{{\cal L}}(\ze ) = {\bf G}_1  {\bf L}_1(\ze ) {\bf G}_1^{-1}
= \left(\begin{array}{cc} \tilde{\rP}(\ze) & \tilde{\rR}(\ze) \\ 
\tilde{\rQ}(\ze) & -\tilde{\rP}(\ze) \end{array}\right)
$$ 
and see that the discrete Lax equation (\ref{dlax})  is transformed to 
\beq\label{ndlax} 
{\cal L}{\cal M}={\cal M}\tilde{{\cal L}}, 
\eeq 
which gives a set of equations equivalent to (\ref{compat}), namely  
\beq\label{bteq} 
\begin{array}{rcl} 
(\ze+v_0) \Big(\tilde{\rP}(\ze) - \rP(\ze)\Big) & = &  -2d_0 \rR(\ze) - \frac{1}{2}\tilde{\rQ}(\ze), \\ 
2d_0\Big( \tilde{\rP}(\ze) + \rP(\ze)\Big) & = & -(\ze+v_0)\rQ(\ze), \\ 
-4d_0\tilde{\rR}(\ze) & = & \rQ(\ze), 
\end{array} 
\eeq 
describing the transformation of the entries of ${\cal L}$. 

By considering the map defined by (\ref{bteq}), we eventually arrive at a proof of Theorem \ref{int}. 

\begin{thm}\label{pmap}  
On the  vanishing level set of the top Casimir, $c_{2\rg+1}=0$, 
the map (\ref{bteq})   is an integrable Poisson map on the space of Lax matrices (\ref{newlax}). 
\end{thm} 
\begin{prf} From the discrete Lax equation (\ref{ndlax}) it is clear that the map given by 
(\ref{bteq}) is isospectral: it leaves the spectral curve (\ref{muspec}) unchanged, and hence 
preserves all the Casimirs (including the constraint $c_{2\rg+1}=0$) and the $\rg$ first integrals 
$H_j$ defined by (\ref{hams}). Thus, for integrability it only remains to show that it is a Poisson map. 
This means that the same bracket relations (\ref{pbs1}), 
(\ref{pbs2}), (\ref{pbs3}), but with tildes,  must hold between the entries of $\tilde{{\cal L}}$, so 
that 
\beq\label{tpbs1} 
\{\, \tilde{\rP}(\ze), \tilde{\rP}(\eta)\,\}  =  0  =  
\{\, 
\tilde{\rR}(\ze), 
\tilde{\rR}(\eta)\,\}, 
\eeq 
and so on. 
All six bracket relations can be checked directly by  substituting for $\tilde{\rP}$,  $\tilde{\rQ}$,  $\tilde{\rR}$ 
in terms of $\rP$, $\rQ$, $\rR$ and then using the brackets between the original 
polynomials (without tildes), but this is extremely tedious, and it is possible 
to bypass most of these calculations. To begin with note that, from (\ref{dvdef}) 
and the 
first bracket in (\ref{pbs2}),  we have 
\beq\label{Pd} 
\{\, d_0,\rP(\eta)\,\}  = \frac{1}{4} \mathrm{res}\, \frac{\{\,\rP(\eta),\rQ(\ze)\,\}}{\ze^{\rg+1}}\Big|_{\ze=0}
=\frac{1}{4} \mathrm{res}\, \frac{\rQ(\ze)-\rQ(\eta)}{(\ze-\eta)\ze^{\rg+1}}\Big|_{\ze=0}=0,
\eeq 
since $$\rQ(\ze)-\rQ(\eta)=-4d_0(\ze-\eta)\Big(\ze^{\rg-1}+O(\ze^{\rg-2})\Big),$$
and similar calculations show that 
\beq\label{QRd} 
\{\, d_0,\rQ(\eta)\,\}  = -4d_0, \qquad 
\{\, d_0,\rR(\eta)\,\}  = -1. 
\eeq 
So, substituting $\tilde{\rR}=-\frac{1}{4}d_0^{-1}\rQ$ from (\ref{bteq}), it follows that 
$$ 
\{ \,\tilde{\rR}(\ze), \tilde{\rR}(\eta)\,\} 
=\frac{1}{16}\Big( -d_0^{-3}\rQ(\ze) \{\,d_0,\rQ(\eta)\,\} -d_0^{-3} \{\,\rQ(\ze),d_0\,\}\rQ(\eta)
+d_0^{-2}\{\,\rQ(\zeta),\rQ(\eta)\,\}\Big) =0 
$$
by (\ref{QRd}) and the first bracket in (\ref{pbs3}), which 
verifies the second relation in (\ref{tpbs1}). Note that from (\ref{bteq}) 
and (\ref{dvdef}) we may write 
$$ 
v_0=-\mathrm{res}\, \frac{\tilde{\rR}(\ze)}{\ze^\rg}\Big\vert_{\ze=0},
$$ 
from which it follows that 
\beq\label{rv}
\{ \,v_0,\tilde{\rR}(\ze)\,\}=0, 
\eeq 
and, making use of (\ref{Pd}), we also have 
\beq\label{prt} 
\{\,\rP (\ze),\tilde{\rR}(\eta)\,\}= -4d_0\{\,\rP(\ze),\rQ(\eta)\,\} = 2\left( \frac{\tilde{\rR}(\ze)  -\tilde{\rR}(\eta)}{\ze-\eta}\right),
\eeq 
which then shows that 
\beq\label{pv}
\{ \,v_0,\rP(\ze)\,\} = \mathrm{res} \, \frac{2(\tilde{\rR}(\ze) - \tilde{\rR}(\eta))}{\eta^\rg (\ze - \eta)}\Big\vert_{\eta=0}=2.
\eeq 
Then from the second and third equations in (\ref{bteq}) we may write 
$$ 
\tilde{\rP}(\ze)=-\rP(\ze)+2(\ze + v_0)\tilde{\rR}(\ze),
$$ 
so that from the first bracket in (\ref{pbs1}), and the fact that all the entries of $\tilde{\rR}$ are in involution, together with (\ref{rv}), 
we have 
$$ \begin{array}{rcl}
\{ \,\tilde{\rP}(\ze), \tilde{\rP}(\eta)\,\} 
& = & -2\Big( \{\,\rP(\ze) , (\eta + v_0)\tilde{\rR}(\eta)\,\} + \{\,(z+v_0)\tilde{\rR}(\ze),\rP(\eta)\,\}\Big)\\ 
& = &  -2\Big( \{\,\rP(\ze) ,  v_0\,\} \tilde{\rR}(\eta) +\{\,v_0,\rP(\eta)\,\} \tilde{\rR}(\ze) + (z+v_0)\{\,\tilde{\rR}(\ze),\rP(\eta)\,\} 
+(\eta+v_0)\{\,\rP(\ze),\tilde{\rR}(\eta)\,\}\Big) \\ 
&=&0,
\end{array}
$$
by (\ref{prt}) and (\ref{pv}), which verifies the first bracket in (\ref{tpbs1}). 
Finally, using the same set of bracket identities, we are able to check that 
\beq \label{prtbr}
\{\,\tilde{\rP} (\ze),\tilde{\rR}(\eta)\,\}= 
-\{\,\rP (\ze),\tilde{\rR}(\eta)\,\}= - 2\left( \frac{\tilde{\rR}(\ze)  -\tilde{\rR}(\eta)}{\ze-\eta}\right),
\eeq 
corresponding to the shifted version of the second bracket relation in (\ref{pbs2}). 
It is not necessary to verify directly that the remaining three bracket relations are preserved by the map, 
since they follow from observing that, given the $\rg$ pairs of coordinates $(\tilde{x}_i,\tilde{y}_i)$
defined by 
$$ 
\tilde{\rR}(\ze) =\prod_{j=1}^\rg (\ze - \tilde{x}_j) , \qquad \tilde{y}_i =\tilde{\rP}(\tilde{x}_i), 
$$ 
for $i=1,\ldots,\rg$, the already verified relations (\ref{tpbs1}) and (\ref{prtbr}) imply that 
$$
\{\, \tilde{x}_i,\tilde{x}_j\,\} = 0 =\{\, \tilde{y}_i,\tilde{y}_j\,\} , \qquad 
\{\, \tilde{y}_i,\tilde{x}_j\,\}=2\delta_{ij}.
$$
Hence the map defined by (\ref{bteq}) restricts to a canonical transformation on each 
symplectic leaf, and it preserves all the Casimirs, so it is a Poisson map. 
This also proves Theorem  \ref{int}, since $\pi_0^{(j)}$ together with  $\rho_{-1}^{(j)}$ (the coefficients of $\rR$) 
for $j=1,\ldots,\rg$ also provide coordinates on each symplectic leaf, so that the map (\ref{phihat}) which is written 
in these coordinates is symplectic and has $\rg$ first integrals in involution, as in (\ref{hams}). Thus $\hat\varphi$ 
is integrable, as is the conjugate map $\varphi$ given by (\ref{phi}).
\end{prf}

\section*{Appendix: Identities for determinants of Hankel type}
There are various ways to derive the classical formulae (\ref{hankform}), 
and the expression for $d_n$ in particular; 
see the proof of Theorem A in \cite{wall}, for instance. However, 
here we present determinantal formulae that yield the latter expression 
directly from the three-term relation (\ref{3term}). 

For convenience, we introduce some notation:  take the column vectors 
${\bf c}_j = (s_j,\ldots,s_{j+n-2})^T$,
${\bf c}_j' =(s_j,\ldots,s_{j+n-1})^T$ of sizes $n-1$ and $n$, respectively, and let   
$$ 
\Delta_n^{**} 
=
\left| \begin{array}{cccccc} 
s_0       & s_1      &  \cdots & s_{n-3} & s_{n-1} &  s_{n}   \\ 
 s_1      & \vdots  &           & \vdots  &\vdots &\vdots      \\
\vdots   & \vdots &          &   \vdots  &\vdots &\vdots  \\        
s_{n-1} &  s_n     & \cdots & s_{2n-4} & s_{2n-2} & s_{2n-1}  
\end{array} 
\right| =|{\bf c}_0'\cdots {\bf c}_{n-3}' {\bf c}_{n-1}'{\bf c}_{n}'|
$$ 
denote a size $n$ determinant of Hankel type  with the column 
${\bf c}_{n-2}'$ omitted.
Upon using  (\ref{qform}) to calculate the coefficient of order $X^{n-1}$ in (\ref{3term}),  it follows 
that $d_n$ is given as a linear combination of four terms, that is  
$$ 
d_n =  \frac{\Delta_{n-1}^{**}}{\Delta_{n-1}} -\frac{\Delta_{n}^{**}}{\Delta_{n}}
-\frac{\Delta_{n-1}^{*}}{\Delta_{n-1}}
\left( \frac{\Delta_{n-1}^{*}}{\Delta_{n-1}}-\frac{\Delta_{n}^{*}}{\Delta_{n}} \right). 
$$ 
First of all, observe that the Desnanot-Jacobi identity, also known as Dodgson condensation \cite{dodgson}, 
yields the formula 
$$ 
\Delta_n\Delta_{n-2} = \Delta_{n-1}\Delta_{n-1}' - (\Delta_{n-1}^*)^2, 
$$ 
where $\Delta_{n-1}'$ is a matrix whose first principal minor of size $n-2$ is $\Delta_{n-2}$, 
namely 
$$ 
\Delta_{n-1}' = 
\left| \begin{array}{ccccc} 
s_0     & s_1 & \cdots & s_{n-3} & s_{n-1}   \\ 
s_1     &         & \iddots &        \vdots  & \vdots      \\
\vdots & \iddots &          &           \vdots  &\vdots   \\        
s_{n-3}     &  \cdots & \cdots & s_{2n-6} & s_{2n-4} \\   
s_{n-1}     &  \cdots & \cdots & s_{2n-4} & s_{2n-2}
\end{array} 
\right|. 
$$
Then the above formula for $d_n$  implies that 
$$ 
d_n - \frac{\Delta_n\Delta_{n-2}}{\Delta_{n-1}^2} =
 \frac{\Delta_n\Delta_{n-1}^{**}-\Delta_{n}^{**}\Delta_{n-1}
-\Delta_{n}\Delta_{n-1}'+\Delta_{n}^{*}\Delta_{n-1}^{*}}{\Delta_{n-1}\Delta_n},
$$  
and this can be simplified further by 
introducing 
$$ 
\Delta_{n-1}'' = |{\bf c}_0\cdots{\bf c}_{n-3}{\bf c}_n|, 
$$ 
and then considering a  determinant of size $2n-1$, namely  
$$ 
\left| \begin{array} {ccccccccc}
 {\bf c}_0' &\cdots &{\bf c}_{n-3}'& {\bf c}_{n-2}' &{\bf c}_{n-1}' & {\bf c}_{n}'& \mathbf{0}& \cdots &  \mathbf{0}\\ 
\mathbf{0}& \cdots &  \mathbf{0} & {\bf c}_{n-2} &{\bf c}_{n-1} & {\bf c}_{n}& {\bf c}_0 &\cdots &{\bf c}_{n-3} 
\end{array} \right|=0 ,
$$ 
which can be seen to vanish from elementary row operations. 
Performing the Laplace expansion of the latter  determinant into products of blocks 
of sizes $n$ and  $n-1$  gives just 
three non-zero terms, producing the identity 
$$ 
\Delta_{n}^{**}\Delta_{n-1} - 
\Delta_{n}^{*}\Delta_{n-1}^{*}+ \Delta_{n}\Delta_{n-1}''=0, 
$$ 
which reduces the expression for $d_n$ to 
$$ 
\Delta_{n-1} \left(d_n - \frac{\Delta_n\Delta_{n-2}}{\Delta_{n-1}^2} \right) = 
\Delta_{n-1}''-\Delta_{n-1}'+ \Delta_{n-1}^{**}. 
$$ 
Finally, to see that the right-hand side above vanishes, shift $n\to n+1$, and then note that 
$$ 
\Delta_n'' - \Delta_n'+\Delta_n^{**} = 
\left| \begin{array}{ccccc} 
{\bf c}_0 & \cdots & {\bf c}_{n-3} & {\bf c}_{n-1} & {\bf c}_n \\ 
s_{n-1} & \cdots & s_{2n-4} & s_{2n-2} & s_{2n-1}  
\end{array}\right|
- 
\left| \begin{array}{cccc} 
{\bf c}_0 & \cdots & {\bf c}_{n-2}  & {\bf c}_n \\ 
s_{n} & \cdots  & s_{2n-2} & s_{2n}  
\end{array}\right|
+ 
\left| \begin{array}{ccc} 
{\bf c}_0 & \cdots  & {\bf c}_{n-1}  \\ 
s_{n+1} & \cdots  & s_{2n}  
\end{array}\right|, 
$$
where the latter are the only three non-zero terms that appear 
in the  sum 
$$ 
\sum_{j=0}^n (-1)^j \left| \begin{array}{lclcl} 
{\bf c}_0 & \cdots & \widehat{{\bf c}}_{j} & \cdots  & {\bf c}_n \\ 
s_{j+1} & \cdots & \widehat{s}_{2j+1} 
& \cdots & s_{n+j+1}  
\end{array}\right|=0
$$ 
(with the hat denoting an omitted column), which 
is seen to be identically zero by expanding about the last row.

\small 
\noindent \textbf{Acknowledgments:} This research was supported by 
Fellowship EP/M004333/1  from the Engineering \& Physical Sciences Research Council, UK. 
I am grateful to the  
School of Mathematics and Statistics, University of New South Wales, for 
hosting me as a Visiting Professorial Fellow with funding from the Distinguished Researcher Visitor Scheme, and 
to John Roberts and Wolfgang Schief, who provided additional support during my 
time in Sydney. I would also like to thank Shihao Li for inviting me to visit the University of Melbourne, 
where we had many enlightening discussions. This work is dedicated to Jon Nimmo, who first 
made me aware of the results in \cite{ch} and suggested that there should be a way to extend them  
to higher order Somos sequences. Jon was an expert on symmetric functions and associated 
determinantal formulae for the solutions of integrable systems \cite{nimmo}, and I like 
to think he would have appreciated the identities in the appendix.


\end{document}